\documentclass{amsart}

\usepackage{amssymb}
\usepackage{amsmath}
\usepackage{graphicx}
\usepackage{float}
\usepackage{hyperref}
\newtheorem{theorem}{Theorem}[section]
\newtheorem{lemma}[theorem]{Lemma}
\newtheorem{prop}[theorem]{Proposition}
\newtheorem{cor}[theorem]{Corollary}

\theoremstyle{definition}
\newtheorem{definition}[theorem]{Definition}
\newtheorem{example}[theorem]{Example}

\theoremstyle{remark}
\newtheorem{remark}[theorem]{Remark}
\newtheorem{notation}[theorem]{Notation}
\numberwithin{equation}{section}

\newcommand{\rs}{\mathrm{res}}
\newcommand{\SE}{\mathcal{S}\hspace{-0.04cm}\mathit{e}}
\newcommand{\ttp}{\tau}
\newcommand{\st}{\mathrm{Stab}}
\newcommand{\lbe}{L_{\mathit{be}}}
\newcommand{\adm}{\mathcal{A}}
\newcommand{\mm}{\mathcal{M}}
\newcommand{\ese}{\mathbb{E}_6}
\newcommand{\ord}{\mathrm{ord}}

\begin{document}

\title{Construction of symmetric cubic surfaces}

\author{Michela Brundu}
\address{Department of Mathematics and Geosciences}
\curraddr{Via Valerio 12/A, 34127 Trieste, Italy}
\email{brundu@units.it}
\thanks{}

\author{Alessandro Logar}
\address{Department of Mathematics and Geosciences}
\curraddr{Via Valerio 12/A, 34127 Trieste, Italy}
\email{logar@units.it}
\thanks{}

\author{Federico Polli}
\address{Area Science Park}
\curraddr{Padriciano 99, 34149, Trieste, Italy}
\email{polli.f@outlook.it}
\thanks{}

\subjclass[2020]{14J50, 14Q10}

\date{}

\dedicatory{}

\begin{abstract}
  We consider the action of the group $\mathrm{PGL}_4(K)$ on the
  smooth cubic surfaces of $\mathbb{P}^3_K$ ($K$ an algebraically closed field of characteristic zero). We classify, in an explicit way, all the smooth cubic surfaces
  with non trivial stabilizer, the corresponding stabilizers and 
  obtain a geometric description of each group in terms of permutations
  of the Eckardt points, of the $27$ lines or of the $45$ tritangent
  planes.
\end{abstract}

\maketitle

\section{Introduction}

\bigskip

Traditionally, as it is claimed in  the Segre's book~\cite{segre},
``the study of the general cubic surface
dates from 1849, in which year the $27$ lines were found by Cayley and
Salmon''. Nevertheless,  as pointed out by~\cite{dol}, cubic surfaces were considered for the first time in
 a work of Pl\"ucker, which dates of 1829. Certainly,
the subject is very old and classic.
A wide historical overview on the theme can be found
in~\cite{dol} (see also Nguyen's  Thesis~\cite{nc1}).

Nevertheless,
the beauty and richness of the properties of these surfaces 
inspired and till inspire new researches on the subject. Consider
that, even in very recent years, an entire issue of the Journal
``Le Matematiche" has been devoted to the cubic surfaces.
It is hard to draw up complete references  on the topic,  we point
out however the ample bibliography  in the  paper~\cite{rs}  
(where an interesting list of open problems is
given) and in the book~\cite{dol}, where an entire chapter is devoted
to give a modern view to the cubic surfaces.

In addition, in the last years several authors have considered the
problem of classifying cubic surfaces over finite fields (see,
for instance,~\cite{dodu},~\cite{bk} and the references given there).

Concerning Segre's investigation, in ~\cite{segre} he described,
in particular, the groups of symmetries of the smooth cubic
surfaces and gave the list of them. He also
realized that non-trivial symmetries are
connected to the existence of \emph {Eckardt points} (also known
in the literature  as \emph {star points}), i.e.\ points
which are the intersection of three coplanar lines of the surface.

In more recent years, many authors studied both these topics. In
the paper~\cite{hosoh}, Hosoh reconsidered the problem of possible
automorphisms of cubic surfaces and, starting from their description
 as the blow-ups of 6 generic points of the plane,
obtained all the groups. In particular, he pointed out some mistakes
in the book of Segre and found a further surface (whose
symmetric group is $C_8$), that was missed in~\cite{segre}. In 2012
Dolghacev, in Chapter~9 of his book~\cite{dol}, gave another, complete
description, of the automorphisms groups and the possible types of the
corresponding surfaces.

Another problem that has received the attention of many authors is the
 determination of the moduli space of cubic surfaces
 (see~\cite{brd}, \cite{dol} Section 9.4.5) and
also~\cite{nc2}, where the description of the moduli space is
connected to the presence of Eckardt points.

In the present paper, we consider one more time the
classification of the automorphisms of smooth
cubic surfaces in the three-dimensional projective space.

We follow an approach which allows us to obtain a four dimensional
family of cubic surfaces which, from one hand, parametrizes
(up to a projectivity) all the smooth cubic surfaces and,
from the other hand, allows us to get the explicit equations of the $27$ lines
for each cubic surface of the family.
{From} the knowledge of the lines then we can easily determine the
$45$ tritangent planes and the possible Eckardt points. Therefore
we stratify the initial four-dimensional family
into subfamilies which parametrize cubic surfaces with all
the possible configurations of Eckardt points.

Consequently,  we study the stabilizer (w.r.t.\ the action of
$\mathrm{PGL}_4(K)$) of every surface of the obtained families,
considering the
corresponding permutations of the Eckardt points, of the $27$ lines or
of the tritangent planes.

The paper is organized as follows: in Section 2 we introduce several
tools in order to get the equation of a suitable family of cubic surfaces which parametrizes
(up to a projectivity) all the smooth cubic surfaces of $\mathbb{P}^3$ and
for which the $27$ lines are explicitly determined. 

As a consequence, in
Section 3 we are able to re-obtain many known results on the configuration
of the possible Eckardt points which allow us to subdivide the surfaces
into several subfamilies. 

The knowledge of the lines 
is the basic point for the contruction given in Section~4 (where we compute
the projectivities of $\mathbb{P}^3$  stabilizing the considered
surfaces) and for the study given in Section~5 (where we complete the determination of
the stabilizers, introducing further subfamilies of surfaces). In this way, we get the list of
 the stabilizers which is clearly the same list
obtained by Segre in~\cite{segre} (with the above mentioned exceptions),
by Hosoh in~\cite{hosoh}, by Dolghacev in~\cite{dol}.
 Our approach (which is very elementary) gives, in addition, a uniform way to understand the
classification and, furthermore, all the object we manipulate (surfaces, Eckardt points,
lines, tritangent planes, projectivities, \dots) are totally explicit.

The final Section, making use of the previous constructions, presents a geometric interpretation of the automorphisms groups
in terms of permutations of lines, tritangent planes and Eckardt points.

While we are aware that  it is certainly not easy to obtain truly
new results after so many years of research on the subject,
as far as we know, the approach we propose is not present in the
literature and provides some new tools to treat the subject.

In order to obtain these results, we have intensively used packages
of symbolic computation (see~\cite{CoCoA-5} and~\cite{sage}) and, in
particular, we have implemented some specific Sage software, available at the repository: \\
\verb+https://github.com/FedericoPolli/Simmetries_Of_Cubic_Surfaces+. \\
The interested reader can download from this site all the procedures we used and also (in the directory \verb+computations+)
several Sage sections which contain the construction of all the possible families of surfaces according to their Eckardt points  and, for each  family, the  elements of the automorphisms group, as used in the present paper.


\section{Preparatory results}
\label{Sec2}
Let us set the following notation:  $K$ is an algebraically closed field of characteristic zero, $\mathbb{P}^3_K$ is the projective space on $K$, whose homogeneous coordinates are $[x,y,z,t]$, and $\mathrm{PGL}_4(K)$ is the projective general linear group ($4 \times 4$ invertible matrices, up to a scalar) acting in the canonical way on $\mathbb{P}^3_K$. We refer to it as the {\em group of projectivities} of $\mathbb{P}^3_K$.

Finally,  let $K[x, y, z, t]$ be the ring of
polynomials in four variables over $K$. 

\medskip

If $F \in K[x, y, z, t]$ is a
homogeneous polynomial of degree three, the set of zeroes of $F$, denoted by $S = V(F)$, is a {\em cubic surface}  in $\mathbb{P}^3_K$. 

It is well known that, if $S$ is smooth,
it contains $27$ lines having a precise configuration
(see~\cite{hart}, chap.~V, 4) that we briefly recall. They can be
labelled by:
\begin{equation}
\label{27rette}
E_i, G_j  \ (\mbox{for } i = 1, \dots, 6), \ \ F_{ij} \
(\mbox{for } 1 \leq i < j \leq 6)
\end{equation}
and they intersect according to the rules:
$E_i$ intersects $G_j$ if and only if $i \not = j$, $E_i$ or $G_i$
intersects $F_{hk}$ if and only if $i \in \{h, k\}$, $F_{ij}$ intersects
$F_{hk}$ if and only if $i, j, h, k$ are all distinct.

Cubic surfaces are parametrized by $\mathbb{P}^{19}_K$, hence a space describing them, up to projectivities, is four-dimensional. In the literature one can find many different ways to introduce it. The construction of the four dimensional family we present here has the advantage to explicitly give   all the lines of the cubic surfaces.

In order to do so, we summarize  the  general approach in  \cite{bl} even if, in the present paper, we restrict ourselves to the smooth case.

\begin{definition} An $L$-set is a quintuple $(l_1, l_2, l_3, l_4, l_5)$
  of lines of $\mathbb{P}^3_K$  such that   $l_2$ intersects $l_1$, $l_3$ and $l_5$, while $l_4$ intersects only $l_1$ and $l_3$ and there are no further intersections.
\end{definition}

\begin{lemma}
  If $L_1 = (l_1, \dots, l_5)$ and $L_2=(l_1', \dots, l_5')$ are two
  $L$-sets, then there exists a unique projectivity of $\mathbb{P}^3_K$
  which sends $l_i$ to $l_i'$, for $i=1, \dots, 5$.
  \label{projLset}
\end{lemma}

\begin{notation} We denote such a projectivity by
  $M(L_1, L_2)\in \mathrm{PGL}_4(K)$.  
\end{notation}

It is showed that every smooth cubic surface contains an $L$-set (and precisely $25,920$).  As a consequence,
if we choose a specific $L$-set  then the family of cubic surfaces
containing it represents all the smooth cubic
surfaces, up to a projectivity.

In the sequel, we will choose a specific  $L$-set as follows.

\begin{definition}
We call  \emph{basic $L$-set}, and denote it by $L_b$, the following quintuple
\begin{equation}
  \label{pettine}
L_b = (l_1, l_2, l_3, l_4, l_5) = \left((y, z),\ (x,y),\ (x,t),\ (x-z,y-z),\ (x-y, z+t) \right).
\end{equation}
\end{definition}

The family of all cubic surfaces passing through  $L_b$ is the four dimensional linear system given by
\begin{equation}
\begin{split}
& a(2x^2y-2xy^2+xz^2-xzt-yt^2+yzt)+b(x-t)(xz+yt)+\\
& c(z+t)(yt-xz)+d(y-z)(xz+yt)+g(x-y)(yt-xz) =  0.
\end{split}
\label{cubicaA}
\end{equation}
The parameters $a,b,c,d,g$ give singular surfaces if and only if satisfy  $\sigma = 0$, where
\begin{equation}
\label{luogoSingA}
\begin{split}
  \sigma =& c(a+b-c)(2 a+b-d)(a-c-d)(a+c+g)\cdot\\
  &(a+c-g)(4 a c-g^{2})(a^{2}+a c-2 a d+a g+d^{2}-d g)\cdot\\
  &(a^{2}+2 a b+a c-a g+b^{2}-b g)\cdot \\
	&(4 a^{2}+3 a b-4 a c-3 a d-b c-2 b d+b g+c d+d g)\cdot \\
	&(4 a^{3}+4 a^{2} b-8 a^{2} c-4 a^{2} d+a b^{2}-4 a b c-2 a b d+2 a b g+4 a c^{2}+\\
	&4 a c d+a d^{2}+2 a d g+b^{2} c+b^{2} g+2 b c d-2 b c g+c d^{2}-2 c d g-d^{2} g)
\end{split}
\end{equation}
and we briefly say that it is the {\em singular locus} inside the family (\ref{cubicaA}).

\medskip

\begin{definition} 
A plane $\pi$ is called \emph{tritangent} to a smooth cubic surface $S$ if  $\pi \cap S$ consists of three  (distinct)  lines. \\
If $r$ and $s$ are two meeting lines of a smooth cubic surface $S$, then the plane containing them  is tritangent  $S$  and the third line  will be denoted by $\rs(r, s)$ and  called the \emph{residue line} of $r$ and $s$.
\end{definition}

If $r$ and $s$ are as above, the tritangent plane  containing them will be denoted  by the triplet of lines 
$(r,s, \rs(r,s))$ or simply by $(r,s)$.

\medskip
The properties of incidence of the $27$ lines (see (\ref{27rette})) allow us to realize that there are $45$ tritangent planes to a smooth cubic surface $S$ (see  Table~\ref{45planes}).

\medskip


\begin{table}[H]
  \caption{The 45 tritangent planes. }
  \label{45planes}
  \begin{center}
    $
\small{
\begin{array}{|l|l|l|l|l|}\hline
\ttp_{1} : (E_1, G_2) & \ttp_{2} : (E_1, G_3) & \ttp_{3} : (E_1, G_4) &
\ttp_{4} : (E_1, G_5) & \ttp_{5} : (E_1, G_6)\\
\ttp_{6} : (E_2, G_1) & \ttp_{7} : (E_2, G_3) & \ttp_{8} : (E_2, G_4) &
\ttp_{9} : (E_2, G_5) & \ttp_{10} : (E_2, G_6)\\
\ttp_{11} : (E_3, G_1) & \ttp_{12} : (E_3, G_2) & \ttp_{13} : (E_3, G_4) &
\ttp_{14} : (E_3, G_5) & \ttp_{15} : (E_3, G_6)\\
\ttp_{16} : (E_4, G_1) & \ttp_{17} : (E_4, G_2) & \ttp_{18} : (E_4, G_3) &
\ttp_{19} : (E_4, G_5) & \ttp_{20} : (E_4, G_6)\\
\ttp_{21} : (E_5, G_1) & \ttp_{22} : (E_5, G_2) & \ttp_{23} : (E_5, G_3) &
\ttp_{24} : (E_5, G_4) & \ttp_{25} : (E_5, G_6) \\
\ttp_{26} : (E_6, G_1) & \ttp_{27} : (E_6, G_2) & \ttp_{28} : (E_6, G_3) &
\ttp_{29} : (E_6, G_4) & \ttp_{30} : (E_6, G_5) \\
\ttp_{31} : (F_{12}, F_{34}) &\ttp_{32} : (F_{12}, F_{35}) & \ttp_{33} : (F_{12}, F_{36}) &
\ttp_{34} : (F_{13}, F_{24}) & \ttp_{35} : (F_{13}, F_{25})\\
\ttp_{36} : (F_{13}, F_{26}) & \ttp_{37} : (F_{14}, F_{23}) & \ttp_{38} : (F_{14}, F_{25}) &
\ttp_{39} : (F_{14}, F_{26}) & \ttp_{40} : (F_{15}, F_{23})\\
\ttp_{41} : (F_{15}, F_{24}) & \ttp_{42} : (F_{15}, F_{26}) & \ttp_{43} : (F_{16}, F_{23}) &
\ttp_{44} : (F_{16}, F_{24}) & \ttp_{45} : (F_{16}, F_{25})\\ \hline
\end{array}
}
$
  \end{center}
\end{table}

\noindent
Here we  list a few simple properties about  the  $27$  lines of a smooth cubic surface $S$.

\begin{remark}
  \label{banale}
  Let $r$, $s$, $t$ be three lines on $S$. If $r$ and
  $s$ are incident and $t$ is skew with $r$ and $s$, then $t$ meets
  $\rs(r, s)$.
\end{remark}

\begin{lemma}
Let $\{s_1,s_2\}$ and $\{r_1,r_2,r_3\}$ two sets of skew lines on a smooth cubic surface
such that $s_i \cap r_j \ne \emptyset$, for all $i,j$. Then the line $t =\rs(\rs(r_1,s_1), \rs(r_2,s_2))$
  intersects $r_3$.
  \label{treRette}
\end{lemma}
\begin{proof}
  The lines $\rs(r_1, s_1)$ and $\rs(r_2, s_2)$ are skew with $r_3$, so,
  by Remark~\ref{banale}, $t$ intersects $r_3$.
\end{proof}

Again from the  incidence properties  (see (\ref{27rette})), a simple argument shows that, given two skew line $s_1$ and $s_2$ on $S$, there are five lines of $S$  intersecting $s_1$ and $s_2$. The next result 
is more precise  and leads to determine by residuality all the lines of $S$.

\begin{prop} 
Let $\{s_1,s_2\}$ and $\{r_1,r_2,r_3,r_4\}$ two sets of skew lines on a smooth cubic surface $S$ 
such that $s_i \cap r_j \ne \emptyset$, for all $i,j$. Then all the $27$ lines of $S$ can be determined by
residuality.
\label{quintaRetta}
\end{prop}
\begin{proof}
The line $t = \rs(\rs(r_1,s_1), \rs(r_2,s_2))$ intersects $r_3$ and $r_4$ from Lemma~\ref{treRette}.
 Consider now the line $u= \rs(\rs(r_3,s_2), \rs(r_4,t))$.
Then, again from Lemma~\ref{treRette}, $u$ intersects $s_1$. Since $s_1$ and $u$ are skew with $s_2$, we get from
Remark~\ref{banale} that the line $r_5 = \rs(s_1, u)$ intersects $s_2$, and clearly also $r_1$. \\
Now it is straightforward to see that $\{s_1,s_2\}$ and $\{r_1,r_2,r_3,r_4, r_5\}$ give rise to the remaining $20$ lines of $S$ in the following way: 
\[
\rs(r_i, s_j) \quad \hbox{for} \quad  i=1, \dots, 5, \quad j = 1, 2
\]
\[
\rs(\rs(r_i,s_1), \rs(r_j,s_2))  \quad \hbox{for} \quad  i = 1, \dots, 4, \quad j = i+1, \dots, 5.
\]
\end{proof}

This result leads us to introduce the notion of another  useful sextuple of lines.

\begin{definition}
  A $6$-uple of lines $(l_1, \dots, l_6)$ of a smooth cubic surface $S$
  such that $L=(l_1, \dots, l_5)$
  is an $L$-set and $l_6$ intersects $l_2$ and $l_4$ and is skew with
  the other lines of  $L$ is called an \emph{extended $L$-set} of $S$.
\end{definition}

\begin{prop} 
Given an extended $L$-set $(l_1, \dots, l_6)$ of a smooth cubic surface $S$,  all the other lines
of $S$  can be determined by residuality, in a unique way.
\label{daExLa27}
\end{prop}
\begin{proof}
Indeed, $l_1, l_3, \rs(l_2, l_5), l_6$ are four skew
lines meeting both $l_2$ and $l_4$. From Proposition~\ref{quintaRetta} we can find all the lines of $S$.
\end{proof}

\begin{cor}
\label{27uniq}
As soon as the labels of the lines of an extended $L$-set are chosen,  the labels of all the 27 lines are uniquely determined.
\end{cor}

This gives immediately the following result.

\begin{cor}
There is a one to one correspondence between extended $L$-sets and the
  permutations  of the $27$ lines preserving the incidence relations.
  \label{corrisp}
\end{cor}

The above facts (Proposition~\ref{daExLa27} and Corollary~\ref{corrisp}) do not hold concerning $L$-sets, as the following result shows.

\begin{prop} 
For any $L$-set $L=(l_1, \dots, l_5)$ of $S$, there exist exactly two lines $l_6$ and $l'_6$ of $S$ such that
 $L_e=(l_1, \dots, l_5, l_6)$ and  $L_e'=(l_1, \dots, l_5, l'_6)$ are two extended $L$-sets.
    \label{2estesi}
\end{prop}
\begin{proof}
There are exactly five (mutually skew) lines meeting both $l_2$ and $l_4$. Clearly, $l_1, l_3, \rs(l_2, l_5)$ are three among them. Setting 
$l_6$ and $l'_6$ the two remaining lines, it follows  that $L_e$ and $L_e'$ are two extended $L$-sets.
\end{proof}

Observe that both $l_6$ and $l'_6$ have the same incidence relations with respect the lines of $L$; in other words, an $L$-set cannot distinguish the two extended $L$-sets which contain it.

\begin{remark}
  Since the number of $L$-sets on a smooth cubic surface is $25,920$, there are $51,840$ extended $L$-sets 
  and therefore, by Corollary~\ref{corrisp}, 
  the group of permutations of the $27$ lines of a cubic surface has
  order $51,840$, in
  accordance to the order of the Weyl group $\ese$
  (see, for instance,~\cite{hart}).
\end{remark}

\begin{notation}    For this reason, from now on, we denote by $\ese$ the {\em group of the permutations of the $27$ lines}. Clearly, each of such permutations preserves the incidence relations among all the lines.
\end{notation}

From Corollary~\ref{27uniq} all the extended $L$-sets are equivalent. 
For this reason, we   assume that the basic $L$-set is
\[
L_b = (l_1, l_2, l_3, l_4, l_5) = (E_1, G_4, E_2, G_3, E_3)
\]
and the \emph{basic extended} $L$-set is
\[
\lbe = (l_1, l_2, l_3, l_4, l_5, l_6) = (E_1, G_4, E_2, G_3, E_3, E_5).
\]
All the smooth cubic surfaces containing the basic $L$-set $L_b$ (whose lines have the equations given 
in (\ref{pettine})) are described
by the equation (\ref{cubicaA}), where $\sigma(a, b, c, d, g)\not = 0$.
However, it is convenient to make the substitution
\[
g = e+f, \quad a = ef/c
\] 
in  (\ref{cubicaA}) and  (\ref{luogoSingA}) and this is possible since $c=0$ gives a singular
cubic surface.  We then obtain the  equation (representing a four-dimensional quadric in $\mathbb{P}^{19}_K$) of the family parametrizing, up to a projectivity, all the smooth cubic surfaces 
of $\mathbb{P}^3_K$ :

\begin{equation}
  \begin{split}
   &b c(t-x)(x z+y t)+c^{2}(z+t)(xz-yt)-cd(y-z)(xz+yt) \\
   &+c(e+ f)(x-y)(xz-yt)-e f(2 x^{2} y-2 x y^{2}+x z^{2}-x z t+y z t-y t^{2})=0
  \end{split}
\label{cubicaB}
\end{equation}
and the corresponding singular locus
\begin{equation}
  \begin{split}
    &\Sigma_0 = c(c-f)(-c+e)(c+f)(c+e)(-e+f)(-cd+cf+ef)(-cd+ce+ef)\cdot\\
    & (-c^{2}-cd+ef)(bc-cf+ef)(bc-ce+ef)(bc -cd+2ef)(bc-c^{2}+ef)\cdot\\
    & (bc^{2}+c^{2}d+bcf-2c^{2}f-cdf+2ef^{2})
      (bc^{2}+c^{2}d+bce-2c^{2}e-cde+2e^{2}f)\cdot\\
    & (-bc^{3}-2bc^{2}d+c^{3}d+bc^{2}e+c^{2}de+bc^{2}f+c^{2}df+
    3bcef-4c^{2}ef\\
    & -3cdef+4e^{2}f^{2}).
  \end{split}
  \label{luogoSingB}
\end{equation}

Following the proof of Proposition~\ref{2estesi}
, we can determine exactly two  lines $l_6$ and $l_6'$ that complete $L_b$ to an extended $L$-set. We choose as $l_6$ (or, equivalently,  $E_5$) the line having Pl\"ucker coordinates:
\begin{equation}
  \begin{split}
    & [0{,}\, (f-c)(cd-cf-ef)(bc-cf+ef){,}\\
    & (c-f)(cd-cf-ef)^2{,}\, (c+f)(bc-cf+ef)^2{,} \\
  &(c+f)(cd-cf-ef)(cf-ef-bc){,}\, 2f(cd-cf-ef)(cf-ef-bc)].
  \end{split}
  \label{E5}
\end{equation}
One can see that $l_6'$ has the Pl\"ucker coordinates obtained from  (\ref{E5}) by exchanging  $e$ and $f$ 
and this is consistent with the fact that (\ref{cubicaB}) and (\ref{luogoSingB})
are symmetric in $e$ and $f$.
After this choice, all the remaining lines can be explicitly obtained from the extended  $L$-set $\lbe$ as in  Proposition~\ref{daExLa27}.

The knowledge of the equations of the $27$ lines for the generic cubic of
the family (\ref{cubicaB}) is the starting point for  the constructions that we are
going to introduce.


\section{Eckardt points and Eckardt families}
\label{Sec3}

\begin{definition}
  A point on a smooth cubic surface $S$ is called an \emph{Eckardt point}
  (see~\cite{segre} or ~\cite{dol}, Section 9.1.4) if it is the intersection of  three (necessarily coplanar) lines of $S$. 
  The  tritangent plane containing these three lines is said  an \emph{Eckardt plane}.
\end{definition}

\begin{notation}
If $P=r \cap s \cap t$ is an Eckardt point of $S$,  it uniquely determines the  Eckardt plane $\pi=(r,s,t)$. Therefore, we shall use the same notation to denote both, i.e.\ we shall write
$P=(r, s, t)$ (or  simply $P=(r, s)$). 
\end{notation}

In this way, Table~\ref{45planes} lists either the tritangent planes and  the {\em possible} Eckardt points (planes).

There are many ways to determine conditions to impose the existence of Eckardt points, but, since we know the equations of all the lines of $S$, as soon as a tritangent plane is given, we can easily obtain the conditions on the coefficients $b, c, d, e, f$ of  (\ref{cubicaB}) forcing the  lines of a certain tritangent plane to meet in a common point.

\begin{example}
Consider the tritangent
plane $\ttp_3 = (E_1, G_4, F_{14})$. It has an Eckardt point if and
only if the three coplanar lines (here given by their Pl\"ucker coordinates) $E_1 = [0{,}\,0{,}\,1{,}\,0{,}\,0{,}\,0]$,
$G_4 = [0{,}\,0{,}\,0{,}\,0{,}\,0{,}\,1]$,
$F_{14}= [0{,}\,bc+c^2+ef{,}\,-c^2-cd+ef{,}\,0{,}\,0{,}\,c(-b+e+f)]$
have a common point. It turns out the  corresponding condition is  
$bc+c^2+ef=0$. Similarly, $\ttp_8 = (E_2, G_4, F_{24})$ has
an Eckardt point if and only if $c^2-cd+ef=0$. 
\end{example}

In the same way, for each tritangent plane $\tau_i$ we can determine a polynomial $P_i \in K[b,c,d,e,f]$
such that $\tau_i$ is an Eckardt plane if and only if  $P_i (b,c,d,e,f)=0$. 

We  collect all the conditions and get a list $\mathcal{Q}$ of the 45 polynomials $P_i(b,c,d,e,f)$, $i=1, \dots, 45$,
here omitted for shortness.

\begin{remark} 
It is not difficult to see that every cubic surface with at least one Eckardt point is projectively equivalent to a cubic surface containing the basic $L$-set (\ref{pettine}) and such that the Eckardt point is $A = \tau_3= (E_1, G_4)$.
The corresponding condition is $\Lambda_1 = \{b = -(c^2+ef)/c\}$ and 
this substitution in equation (\ref{cubicaB}) yields to
a family, say $\SE_1$, representing all the  cubic
surfaces with at least one Eckardt point up to a projectivity
and whose singular locus is given by
$\Sigma_1$, obtained from   (\ref{luogoSingB}) together with the condition
$\Lambda_1$. \\
The above substitution annihilates precisely one of the $45$ polynomials
of $\mathcal{Q}$ (clearly the polynomial $P_3$), hence the general element of $\SE_1$ is a smooth cubic surface
with exactly one Eckardt point.
\label{osservase1}
\end{remark}

The same kind of argument can be used to detect families of surfaces with a larger number of Eckardt points.

\begin{remark} 
A cubic surface $S$ containing two Eckardt points (on a
line of the surface) is projectively equivalent to a surface containing the basic $L$-set $L_b$ (\ref{pettine}) and such that the Eckardt points are
$A = E_1 \cap G_4$ and $D = E_2 \cap G_4$. A direct computation shows that this
family, denoted by $\SE_2$ (see Table~\ref{SEi} and Table~\ref{listofequ}), and its singular locus $\Sigma_2$, can be obtained by the substitutions
\[
\Lambda _2 = \{b = -(c^2+ef)/c, \quad d = (c^2+ef)/c\}
\]
into (\ref{cubicaB}) and (\ref{luogoSingB}). \\
Since, in the  above list $\mathcal{Q}$, no other polynomial but $P_3$ and $P_8$ vanishes,
 we obtain that the general element of $\SE_2$ is a 
smooth cubic surface with exactly two Eckardt points
(contained in a line of the surface). 
\label{osservase2}
\end{remark}

\begin{remark} 
Suppose now that a smooth cubic surface has two Eckardt points  not
contained in one of its lines. Again, we can assume, up to a projectivity, that
the surface passes through the $L$-set $L_b$ and that   $A= E_1 \cap G_4 = \ttp_3$
and $C = E_2 \cap G_3= \ttp_7$ are Eckardt points. The corresponding set of conditions is
\[
\Lambda_3 = \{b = -(c^2+ef)/c, d = (3ef-c^2+cf+ce)/(2c) \}.
\]
When we make this substitution in $\mathcal Q$, we obtain that vanish not only $P_3$ and $P_7$, but also
$P_{34}$. This means that $\tau_{34}= (F_{13}, F_{24}, F_{56})$ is an Eckardt  point, that turns out to be collinear with $A$
and $C$. \\
 Since no other polynomial of
$\mathcal{Q}$ vanishes, we get a family $\SE_3$  which
contains (up to a projectivity) all the smooth cubic surfaces with three
collinear Eckardt points. 
\label{osservase3}
\end{remark}

Collecting the three Remarks above, we have the following result (see also \cite{segre} Ch.~IV, Sect.~XIV and ~\cite{dol}, Proposition 9.1.26).

\begin{prop} The following facts hold.
  \begin{enumerate}
    \item There exist cubic surfaces with precisely one Eckardt point.
    \item If a cubic surface has two Eckardt points, then the line joining them either is contained in the surface (and in this case it cannot contain other Eckardt points) or is not a line of the surface (and in this case it intersects the surface in a third point which is another Eckardt point).
    \item There exist cubic surfaces with precisely two Eckardt points.
    \item There exist cubic surfaces with precisely three Eckard points.
      In this case the three points are collinear.
\end{enumerate}
\label{primieck}
\end{prop}
\begin{proof}
  (2) We have only to
  see that it is not possible to have three collinear Eckardt points contained
  in a line of the surface. This is a consequence of the fact that if
  we impose  that  $(E_1 , G_4)$, $(E_2, G_4)$
  and $(E_3, G_4)$
  are Eckardt points, we obtain that the corresponding cubic surfaces are singular. 
  
  \noindent
 (4)  Suppose that a cubic surface has precisely three Eckardt points
  $A_1$, $A_2$ and $A_3$ which are not collinear. Then, from part (2), the lines $r_1 = A_1+A_2$, $r_2 = A_1+A_3$ and $r_3 = A_2+A_3$ must
  be contained in
  the surface, so the tritangent plane $A_1+A_2+A_3$   contains $r_1,r_2,r_3$ and also the lines
  $\rs(r_1, r_2), \rs(r_1, r_3), \rs(r_2, r_3)$ and this is impossible.
\end{proof}

\begin{cor}
  If a smooth cubic surface contains at least three
  Eckardt points then, up to a projectivity, it is a cubic surface of the
  family $\SE_3$.
\end{cor}
\begin{proof} If the three points  $A_1$, $A_2$ and $A_3$ are  collinear, the situation is described in  Remarks~\ref{osservase2} and \ref{osservase3}: in this case the surface is projectively equivalent to an element of $\SE_3$.
Otherwise, with the same argument of Proposition~\ref{primieck} - (4), at least one of the three lines connecting  $A_1$, $A_2$ and $A_3$
 cannot be contained in the surface. Therefore  it contains a further Eckardt point and, so, the surface contains three collinear Eckardt points. 
\end{proof}

The two results above show that, in order to describe all the possible configurations of $n$ Eckardt points of a smooth cubic surface (with $n \ge 3$), we have  to study the  subfamilies of $\SE_3$.

Therefore, from now on in this Section,  we consider the family $\SE_3$ given by (\ref{luogoSingB}) under the conditions $\Lambda_3$. 
In this setting,  the list $\mathcal{Q}$ specializes in the following way:  clearly $P_3$, $P_7$ and $P_{34}$ vanish 
since they correspond to the Eckardt planes $\ttp_3 = (E_1, G_4, F_{14})$,
$\ttp_7 = (E_2, G_3, F_{23})$ and
$\ttp_{34} = (F_{13}, F_{24}, F_{56})$, as described in~Remark~\ref{osservase3}. Moreover, from Proposition~\ref{primieck} - (2), to impose a fourth Eckardt point could force the surface to have some further Eckardt point. In this case, several polynomials of $\mathcal{Q}$  coincide. The remaining distinct and non zero polynomials are~$14$.  After renaming them as $Q_1, \dots, Q_{14}$, we obtain:
\[
\begin{array}{lll} \smallskip
   Q_1 \ = \ 5c^2-ce-cf+ef; && Q_2 \ = \ 3c^2+ce+cf-ef;\\ \smallskip
   Q_3 \ = \ c^2 + 3ce - cf + ef;  && Q_4 \ = \ c^2 - ce + 3cf + ef;  \\ \smallskip
   Q_5 \ = \ 5c^2 + ce + cf + ef;  && Q_6 \ = \ 3c^2 - ce - cf - ef;  \\ \smallskip
   Q_7 \ = \ c^2 + ce - 3cf + ef;  && Q_8 \ = \ c^2 - 3ce + cf + ef;  \\ \smallskip
   Q_9 \ = \ c^2 + ef;  && Q_{10} \ = \ 3c^2 + e^2;  \\ \smallskip
   Q_{11} \ = \ 3c^2 + f^2;  && Q_{12} \ = \ 2c - e + f;  \\ \smallskip
   Q_{13} \ = \ 2c + e - f;   && Q_{14} \ = \ e + f     
\end{array}
\]
and they are associated to the Eckardt points 
accordingly to Table~\ref{condEck}. 


\begin{table}[H]
  \caption{List of Eckardt planes associated to
   $Q_i=0$.}
  \label{condEck}
  \begin{center}
    $
\begin{array}{|cl|cl|} \hline
  Q_1: &\ttp_1, \ttp_{11}, \ttp_{18} & Q_2: &\ttp_2  \\ 
  Q_3: &\ttp_4, \ttp_{28}, \ttp_{35} &
           Q_4: &\ttp_{5}, \ttp_{23}, \ttp_{36} \\ 
  Q_5:&  \ttp_{6}, \ttp_{13}, \ttp_{17} & Q_6: &\ttp_8 \\ 
  Q_7:&  \ttp_9, \ttp_{29}, \ttp_{41} & Q_8:&  \ttp_{10},
  \ttp_{24}, \ttp_{44} \\ 
  Q_9:&  \ttp_{12}, \ttp_{16}, \ttp_{31} & Q_{12}:&  \ttp_{25},
  \ttp_{39}, \ttp_{43} \\ 
  Q_{13}:&  \ttp_{30},  \ttp_{38},  \ttp_{40} & Q_{14}:&
  \ttp_{37}  \\ \hline
  Q_{10}: & \multicolumn{3}{l|}{\ttp_{14}, \ttp_{20}, \ttp_{22}, \ttp_{26},
    \ttp_{33}, \ttp_{42}} \\
  Q_{11}: & \multicolumn{3}{l|}{\ttp_{15}, \ttp_{19}, \ttp_{21}, \ttp_{27},
  \ttp_{32}, \ttp_{45}}\\ \hline
\end{array}
$
\end{center}
\end{table}

\begin{remark}
\label{scambioef}
 The polynomials $Q_4$, $Q_8$, $Q_{11}$ and $Q_{13}$ can be ignored
since they are obtained, respectively, from $Q_3$, $Q_7$, $Q_{10}$ and $Q_{12}$ exchanging
$e$ and $f$ (and we already observed that (\ref{cubicaB}) is invariant
with respect to this exchange).
\end{remark}

We can see, from Table~\ref{condEck}, that there are cubic surfaces with $4$ Eckardt points
(when $Q_2$ or $Q_6$ or $Q_{14}$ are zero) or cubic surfaces with $6$
Eckardt points (when one of the polynomial $Q_1$, $Q_3$, $Q_4$, $Q_5$,
$Q_7$, $Q_8$, $Q_9$, $Q_{12}$, $Q_{13}$ is zero) or cubic surfaces
with $9$ Eckardt points (when $Q_{10}$ is zero).

\begin{lemma}
Let $\mathcal{T}_0$, $\mathcal{T}_1$, $\mathcal{T}_2$  be the  three subfamilies of   $\SE_3$  given by the   conditions  $Q_2 =0$,  $Q_6 =0$, $Q_{14} =0$, respectively. \\
Then,  for every smooth surface $T\in \mathcal{T}_1$ (respectively, $\mathcal{T}_2$) there exists   $S \in \mathcal{T}_0$ such that $S$ and $T$ are projectively
  equivalent, and conversely. 
  \label{3subfam}
\end{lemma}
\begin{proof}
  The cubic surfaces of $\SE_3$ have   $\tau_3 = (E_1,G_4)$, $\tau_7 = (E_2,G_3)$,
$\tau_{34} = (F_{13}, F_{24})$ as  Eckardt points. Moreover, 
 those of $\mathcal{T}_0$ also $\tau_2 = (E_1,G_3)$ and  those of $\mathcal{T}_1$ also $\tau_8 = (E_2,G_4)$. \\
Therefore, on one hand,  $S \in \mathcal{T}_0$  if and only if it has the following Eckardt points
\[
  (l_1, l_2), \ (l_1, l_4),\  (l_3, l_4), \ (\rs(l_1, l_4), \rs(l_2, l_3))
  \]
with respect to the  basic $L$-set
\[
L_b=  (l_1, l_2, l_3, l_4, l_5) = (E_1, G_4, E_2, G_3, E_3).
  \]
On the other  hand, $T \in \mathcal{T}_1$  if and only if it has the following Eckardt points
  \[
  (r_1, r_2),\ (r_1, r_4),\ (r_3, r_4),\ (\rs(r_1, r_4), \rs(r_2, r_3)),
  \]
  with respect to the   $L$-set
 \[
 \overline L =    (r_1, r_2, r_3, r_4, r_5) = (G_4, E_1, G_3, E_2, G_2).
  \]
So, if we put $M=M(\overline L, L_b)$, then for each $T \in \mathcal{T}_1$ the cubic surface $M(T)$ belongs to the family
  $\mathcal{T}_0$. The converse holds by exchanging the roles of $\mathcal{T}_0$ and $\mathcal{T}_1$. \\
  The same argument runs concerning $\mathcal{T}_2$.
\end{proof}

Since the above three subfamilies of $\SE_3$ are projectively equivalent, we can choose one of them as the space paramentrizing cubic surfaces having four Eckardt points. So, we set $\SE_4 = \mathcal{T}_0$.

The conditions and the equation defining $\SE_4$ are contained in Table~\ref{SEi} and in 
Table~\ref{listofequ} (the denominator
$c+e$ can be assumed non zero, since it appears as a factor in the polynomial (\ref{luogoSingB}) defining  the singular locus).

\medskip

Concerning cubic surfaces with $6$ Eckardt points, one can see from Table~\ref{condEck}
that there are $9$ families of this type. With an argument quite similar to the 
proof of Lemma~\ref{3subfam}, we can show that these $9$ families are all projectively equivalent and, so,
 all the  cubic surfaces with $6$ Eckardt
points are parametrized by the subfamily of $\SE_3$ defined, for instance, by the further condition $Q_5=0$. We denote such a family by $\SE_6$.

\medskip
Finally, let us consider the case of $9$ Eckardt points. From Table~\ref{condEck} and 
Remark~\ref{scambioef}, it is enough to consider the only condition
$Q_{10}=3c^2+e^2=0$. This
 gives rise to two  families $\mathcal{T}_1$
and $\mathcal{T}_2$ of cubic surfaces:
\[
\mathcal{T}_1 : \Lambda_3 \cup \{e=\sqrt{-3}c\}, \quad
\mathcal{T}_2 : \Lambda_3 \cup \{e=-\sqrt{-3}c\}
\]
\begin{lemma}
  For every $T_1(c,f)$ in $\mathcal{T}_1$ there exist $c', f'$
  such that the corresponding surface $T_2(c', f') \in \mathcal{T}_2$
  is projectively equivalent to $T_1$ and conversely.
\end{lemma}
\begin{proof} Consider   the $L$-set of $T_1$ given by
$L=(F_{46}, G_6, F_{26}, F_{15}, E_3)$ and  the unique matrix $M = M(L,L_b)$
 (see Lemma~\ref{projLset}).
It can be easily checked that
$M^{-1}(T_1) = T_2(f-c, -c+2 \sqrt{-3} c -f)$.
\end{proof}

The families $\SE_4 = \SE_4(c,e)$, $\SE_6=\SE_6(c, e)$
and $\SE_9=\SE_9(c,f)$ depend on
two parameters, so are one dimensional families.
We specialize
the polynomials of Table~\ref{condEck} with conditions $\Lambda_4$,
$\Lambda_6$ and $\Lambda_9$ respectively, 
in order to see if these
families contain surfaces with more Eckardt points.
We obtain several cubic surfaces with $10$ Eckardt points and
several others with $18$ Eckardt points. Similar arguments to
those discussed above show that, up to projectivities, there is
only one cubic surface $\SE_{10}$ with $10$ Eckardt points obtained
specializing~(\ref{cubicaB}) with the conditions
$\Lambda_{10} = \Lambda_4 \cup \{e=(2-\sqrt{5})c\}$ (whose singular
locus is $\Sigma_{10} = \emptyset$, since $\SE_{10}$ does not depend
on parameters) and only
one cubic surface $\SE_{18}$ with $18$ Eckardt points obtained by
the conditions $\Lambda_{18} = \Lambda_4 \cup \{e = \sqrt{-3}c \}$
(again, the singular locus is $\Sigma_{18} = \emptyset$).

We conclude this section summarizing the above results (see also \cite{segre}, Ch.~IV, Sect.~xiv).

\begin{theorem}
  A smooth cubic surface can have only $0$, $1$, $2$, $3$, $4$,
  $6$, $9$, $10$ or $18$ Eckardt points. If $S$ is a smooth cubic
  surface with $n$
  Eckardt points, then it is projectively equivalent to a cubic surface
  of the family $\SE_n$ (here, for completeness, by $\SE_0$ we denote
  the family~(\ref{cubicaB})) whose singular locus is given by
  $\Sigma_n$.
\end{theorem}

\begin{definition}
The families of cubic surfaces $\SE_n$, where $n = 0,1,2,3,4,6,9,10,18$, will be called {\em Eckardt families}.
\end{definition}

In   Table~\ref{SEi} we list the conditions to impose to (\ref{cubicaB}) in order to obtain the families $\SE_i$. In  the third column, we list the Eckardt points,  labelled accordingly to  Table~\ref{45planes}. In the last column we report the dimension of the Eckardt families.
The families $\SE_1', \SE_1''$ and $\SE_9'$    will be introduced in Section~\ref{Sec5}.


\begin{table}[H]
  \caption{Eckardt families and Eckardt points}
  \label{SEi}
  \small{
  \begin{tabular}{|c|l|l|c|} \hline
    family & conditions & Eckardt points & dim \\ \hline
    $\SE_1$ & $\Lambda_1 : \{b = -(c^2+ef)/c\}$& $\begin{array}{l}\ttp_3\end{array}$ & 3 \\
     $\SE_1'$ & $\Lambda_1 + (\ref{condSe1'})$ & $\begin{array}{l}''\end{array}$ & 1 \\
     $\SE_1''$ & $\Lambda_1 + (\ref{condSe1'})+ (\ref{condSe1''})$ & $\begin{array}{l}''\end{array}$ & 0 \\ \hline
      $\SE_2$ & $\Lambda_2 : \Lambda_1 \cup \{d=(c^2+ef)/c \}$&
      $\begin{array}{l}\ttp_3, \ttp_8 \end{array}$ & 2\\
    $\SE_3$ & $\Lambda_3 : \Lambda_1 \cup \{ d= (3ef-c^2+cf+ce)/(2c)\}$&
                 $\begin{array}{l}\ttp_{3}, \ttp_{7}, \ttp_{34} \end{array}$ & 2 \\
    $\SE_4$ & $\Lambda_4 : \Lambda_3 \cup \{f=c(3c-e)/(c+e)\}$ &
      $\begin{array}{l}\ttp_{3}, \ttp_{7}, \ttp_{8}, \ttp_{34}
        \end{array}$ & 1 \\ \hline
    $\SE_6$ & $\Lambda_6 : \Lambda_3 \cup \{f = -c(5c+e)/(c+e)\}$ &
        $\begin{array}{l}
      \ttp_{3}, \ttp_{6}, \ttp_{7}, \\
      \ttp_{13}, \ttp_{17}, \ttp_{34}
      \end{array}
      $ & 1\\ \hline
    $\SE_9$ & $\Lambda_9 : \Lambda_3 \cup \{e = \sqrt{-3}c  \}$ &
    $\begin{array}{l}
      \ttp_{3}, \ttp_{7}, \ttp_{14}, \\
      \ttp_{20}, \ttp_{22}, \ttp_{26}, \\
      \ttp_{33}, \ttp_{34}, \ttp_{42}
    \end{array}
    $ & 1\\ 
    $\SE_9'$ & $\Lambda_9 + (\ref{condSe9'})$ &
            $\begin{array}{lll}
            & & \\
            &''&\\
            & & 
            \end{array}
      $& 0 \\ \hline
    $\SE_{10}$ & $\Lambda_{10} : \Lambda_4 \cup \{ e = (2-\sqrt{5})c\}$ &
    $\begin{array}{l}
      \ttp_{1}, \ttp_{3}, \ttp_{7}, \ttp_{8}, \\
      \ttp_{11}, \ttp_{12},\ttp_{16}, \\
      \ttp_{18}, \ttp_{31}, \ttp_{34}
     \end{array}
    $ & 0 \\ \hline
    $\SE_{18}$ & $\Lambda_{18} : \Lambda_4 \cup \{e = \sqrt{-3}c \}$ &
    $\begin{array}{l}
      \ttp_{2}, \ttp_{3}, \ttp_{7}, \ttp_{8}, \\
      \ttp_{14}, \ttp_{15}, \ttp_{19}, \ttp_{20},\\
      \ttp_{21}, \ttp_{22}, \ttp_{26}, \ttp_{27}, \\
      \ttp_{32}, \ttp_{33}, \ttp_{34}, \ttp_{37}, \\
      \ttp_{42}, \ttp_{45}
    \end{array}  
    $& 0 \\ \hline
  \end{tabular}
  }
\end{table}

\newpage


\begin{table}[H]
  \caption{Equations of the Eckardt families}
  \label{listofequ}
  \Small{
  \begin{tabular}{|c|l|} \hline
    family & equation \\ \hline
$\SE_0$ & 
$\begin{array}{l}
\\ (\ref{cubicaB}) \\ {}
\end{array}$ 
\\  \hline
$\SE_1$ & 
$\begin{array}{l}
(-x^{2} z  - x z^{2}  - x y t  + y z t  + 2 y t^{2}) c^{2} +  \\
\quad (x y z  - x z^{2} + y^{2} t  - y z t) c d - (x^{2} z  - x y z  - x y t  + y^{2} t )c (e+f) +\\
\quad +  (2 x^{2} y - 2 x y^{2} - x^{2} z  + x z^{2}  - x y t  + y z t) e f 
\end{array}$ 
\\  \hline
$\SE_1'$ & 
$\begin{array}{l}
(x^{2} z  - (w+1) x y z  - x z^{2}  - (w+2) x y t  + 2 y^{2} t  + w y z t  +(w+1)y t^{2}) c^{2} + \\
\quad (-(w+3) x^{2} y + (w+3) x y^{2} + 2 x^{2} z  + (w-2) x y z  + \\
\quad x z^{2} + 2 (w+1) x y t  - (w+2) y^{2} t - wy z t  - (w+1) y t^{2} )c f + \\
\quad ( (w-1) x^{2} y  - (w-1) x y^{2}  + x^{2} z  - x y z  - w x y t + w y^{2} t )f^{2}\\
\mbox{where $w=\sqrt{-1}-1$}
\end{array}$ 
\\  \hline
$\SE_1''$ & 
$\begin{array}{l}
18x^2y - 18xy^2 - 6(w+1)x^2z - 3(w^2 - 4w - 1)xyz + (w^3 - 5w + 6)xz^2+\\
\quad (w^3 - 3w^2 + 7w - 9)xyt - (w^3 + w - 6)y^2t + 3(w^2 - 2w + 1)yzt -\\
\quad (w^3 - 3w^2 + w + 3)yt^2 \\
\mbox{where $w=\sqrt{2}-\sqrt{-1}$}
\end{array}$ 
\\  \hline
$\SE_2$ & 
$\begin{array}{l}
(-x^{2} z + x y z  - 2 x z^{2}  - x y t  + y^{2} t + 2 y t^{2}) c^{2} + \\
\quad(-x^{2} z + x y z  + x y t  - y^{2} t) c (e + f) +\\
\quad  + (2 x^{2} y  - 2 x y^{2} - x^{2} z  + x y z  - x y t  + y^{2} t )e f
\end{array}$
\\  \hline
$\SE_3$ & 
$\begin{array}{l}
(2 x^{2} z + x y z  + x z^{2}  + 2 x y t  + y^{2} t  - 3 y z t  - 4 y t^{2} )c^{2} +\\
\quad ( 2 x^{2} z  - 3 x y z  + x z^{2}  - 2 x y t  + y^{2} t  + y z t) c (e+f) + \\
\quad ( - 4 x^{2} y  + 4 x y^{2}  + 2 x^{2} z  - 3 x y z  + x z^{2}  + 2 x y t  - 3 y^{2} t  + y z t )e f 
\end{array}$
\\  \hline
$\SE_4$ & 
$\begin{array}{l}
(-2 x^{2} z  + 2 x y z - x z^{2}  + x y t  - y^{2} t  + y t^{2}) c^{2} + \\ 
\quad (3 x^{2} y  - 3 x y^{2} - 2 x^{2} z + 2 x y z  - x z^{2}  - 2 x y t  + 2 y^{2} t  + y t^{2}) c e +\\
\quad (- x^{2} y  + x y^{2} + x y t - y^{2} t) e^{2}
\end{array}$
\\  \hline
$\SE_6$ & 
$\begin{array}{l}
(2 x^{2} z  - 4 x y z  + x z^{2}  - 3 x y t + y^{2} t  + 2 y z t  + y t^{2}) c^{2} + \\
\quad (- 5 x^{2} y  + 5 x y^{2} + 2 x^{2} z  - 4 x y z  + x z^{2}  + 2 x y t  - 4 y^{2} t  + 2 y z t  + y t^{2}) c e +\\
\quad (- x^{2} y + x y^{2}  + x y t  - y^{2} t ) e^{2}
\end{array}$
\\  \hline
$\SE_9$ & 
$\begin{array}{l}
(2 x^{2} z  - (w+1) x y z  + x z^{2}  - w x y t  + y^{2} t  + (w-1) y z t  + (w-2) y t^{2} )c + \\
\quad ((w+2) xy (y-x)   
+ 2 x^{2} z - 3 x y z  + x z^{2}  + w x y t  -(w+1) y^{2} t  + y z t) f \\
\mbox{where $w = \sqrt{-3} +1$}
\end{array}$ \\  
\hline
$\SE_9'$ & 
$\begin{array}{l}
  312 \,x^2y - 312 \, xy^2 + 2(5w^3 - 20w^2 + 8w - 32) x^2z - \\
 \quad (w^3 - 56w^2 + 64w - 152)xyz + (5w^3 - 20w^2 + 8w - 32)xz^2 +  \\
\quad 4(w^3 + 9w^2 - 14w - 48)xyt +  (5w^3 - 20w^2 + 8w + 280)y^2t - \\
\quad (9w^3 + 16w^2 - 48w + 88)yzt -    2 (7w^3 - 2w^2 - 20w + 28)yt^2 \\
\mbox{where $w = \sqrt{-1} - \sqrt 3$}
\end{array}$ 
\\  \hline
$\SE_{10}$ & 
$x^{2} y - x y^{2} + 2 x^{2} z - 2 x y z + x z^{2} - 2 x y t + 2 y^{2} t - y t^{2}$\\  \hline
$\SE_{18}$ & 
$3 x^{2} y - 3 x y^{2} - 2 x^{2} z + 2 x y z - x z^{2} - 2 x y t + 2 y^{2} t + y t^{2}$\\  \hline  \end{tabular}
  }
\end{table}


\section{Stabilizers of cubic surfaces of $\SE_n$}
\label{Sec4}

By $\st(S)$ we denote the {\em stabilizer} of a smooth cubic surface $S$ with respect
to the action of the group of  projectivities on $\mathbb{P}_K^3$, i.e. it is the subgroup of $\mathrm{PGL}_4(K)$ defined by
\[
\st(S) =
\{ M \in \mathrm{PGL}_4(K) \mid M(S) = S\}.
\] 
Obviously,  the equality $M(S)= S$ means that the polynomial defining $M(S)$ is a multiple of that defining $S$.

This section is devoted to describe $\st(S)$ as far as $S \in \SE_n$, for all possible $n$. For sake of shortness, 
$\st(S)$ for a generic $S\in \SE_n$ will be denoted by $\st(\SE_n)$.

Assume now that $S$ is a cubic surface containing the basic $L$-set $L_b$
(see~(\ref{pettine})) and take $M \in \st(S)$. Then $M$ sends 
$L_b$ to another $L$-set of~$S$. Hence, in order to find $\st(S)$, it is enough to compute the matrices $M = M(L_b, L)$
where $L$ varies
in the set of $25,920$ $L$-sets and check whether $M(S) = S$. 

The
computation for
the generic cubic~(\ref{cubicaB}), although feasable,  is lengthy and complex, but if we
specialize the
parameters $b, c, d, e, f$ to random numeric values, it becomes
reasonably fast 
and suffices to prove that the stabilizer of the general smooth cubic surface
is trivial.

It is then interesting to detect the cases in which the stabilizer is not trivial. 

\medskip

If $l \subset S$ is a line and  $M \in \st(S)$, then  $M(l)$ is  a line of $S$, so $M$ induces a permutation of the $27$ lines (preserving the incidence relations), naturally defined as
\[
\pi_M = 
\begin{pmatrix}
l_1 & l_2 & \dots & l_{27} \\
M(l_1) & M(l_2) & \dots & M(l_{27}) \\
\end{pmatrix} \in \ese.
\]
In particular, if $L = (l_1, \dots, l_5)$ is an $L$-set, with $\pi_M (L)$ we denote the quintuple $(\pi_M(l_1), \dots, \pi_M(l_5))$ that is clearly still an $L$-set.

\begin{prop} 
If $S$ is any smooth cubic surface, then there is a natural group monomorphism
\[
\phi: \st(S) \longrightarrow \ese \quad \mbox{defined by
      $M \mapsto \pi_M$}.
\]
Consequentely,  $|\st(S)|$ divides $2^7 \cdot 3^4 \cdot 5$.
\label{staine6}
\end{prop} 

\begin {proof}
Obviously,  $\pi_{MN} = \pi_M \circ \pi_N$, for all $M$ and $N$  in $\st(S)$.
 Moreover,  $\phi$ is  injective.  Namely, if $\pi_M$ is the identity permutation, then, in particular, $M(L_b) =L_b$ and therefore $M$ is the identity matrix   by
 Lemma~\ref{projLset}. \\
 The last claim follows from $|\ese| = 51,840 = 2^7 \cdot 3^4 \cdot 5$.
 \end {proof}

From the above result and   Lemma~\ref{projLset}, we have immediately the following fact.

\begin{cor}
For all $M \in \st(S)$ it holds
$M = M(L_b, \pi_M(L_b))$.
\label{inversapi}
\end{cor}

\medskip
The following result describes a relationship between Eckardt points and symmetries of 
cubic surfaces.

\begin{lemma} Let $S$ be a smooth cubic surface and $M \in \st(S)$. The following facts hold:
  \begin{enumerate}
  \item if $\ord(M)=2$, then $S$ has at least one Eckardt
    point;
  \item if $\ord(M) =3$, then $S$ has at least three Eckardt
    points;
  \item if $\ord(M) =5$, then $S$ has at least ten Eckardt
    points.
  \end{enumerate}
  \label{stabConEck}
\end{lemma}
\begin{proof}
 (1)   Since the number of lines of $S$ is odd,
  at least one of them (say $l_2$) is fixed by $M$. As recalled in Section~\ref{Sec2}, there are $5$ tritangent planes
  containing~$l_2$: one is certainly fixed, the other $4$
  are either exchanged in couples
  or fixed. We cover all the possible cases by considering three
  tritangent planes $\pi$,
  $\alpha$ and $\beta$ such that $M(\pi) = \pi$, $M(\alpha) = \alpha$
  and $M(\beta) = \beta$
  or $M(\pi) = \pi$, $M(\alpha) = \beta$ and $M(\beta) = \alpha$.
  Consider now the first
  possibility. Let $\{l_1, \rs(l_1, l_2)\}$,  $\{l_3, \rs(l_2, l_3)\}$,
  $\{l_5, \rs(l_2, l_5)\}$ be the lines of $S \cap \pi$, $S \cap \alpha$
  and $S\cap \beta$
  respectively (different from $l_2$). There are some cases to consider,
  since $M$ can
  fix the lines of these sets or exchange them. In each case, consider
  the skew lines
  $l_1, l_3, l_5$ and the line $l_2$ which intersects them. We can
  complete these four lines
  with a line $l_4$ such that $(l_1, l_2, l_3, l_4, l_5)$ is an $L$-set.
  So, up to
  projectivities, we can assume that this $L$-set is $L_b$. (In particular, we can assume that $S$
  has an equation  given by~(\ref{cubicaB})). If $M$ fixes the three couples
  of lines above, then $M$ sends $L_b$ to either $L_b$ itself or to
  $(E_1, G_4, E_2, F_{12}, E_3)$. If $M$ permutes only the two lines of
  $\beta$ ($E_3$ and
  $F_{34}$) and fixes the lines of $\alpha$, then $M$ sends $L_b$ to either
  $(E_1, G_4, E_2, G_5, F_{34})$ or $(E_1, G_4, E_2, G_6, F_{34})$, and
  so on. In this way we
  collect all the  matrices $M_1, \dots, M_{12}$ which could
  stabilize $S$. For
  $i=1, \dots, 12$, we impose $M_i(S) = S$, obtaining the corresponding
  set of conditions on the parameters of $S$. It turns out that the
  corresponding
  surfaces, when not singular, have at least one Eckardt point. \\
(2) If $M$ has order $3$, we consider
  three cases: either
  one of the lines of $S$ is fixed or there are three lines $r_1, r_2$ and $r_3$
  such that $r_1 \mapsto r_2 \mapsto r_3 \mapsto r_1$ and are coplanar or are skew. Again,
  in the first case we consider the five tritangent planes passing through the fixed lines.
  At least two of them are fixed, so the lines on these planes are also fixed. In particular,
  we can find three lines $l_1$, $l_2$, $l_3$ on $S$ fixed by $M$ and such that $l_1$ and
  $l_3$ are skew and $l_2$ meets $l_1$ and $l_3$. Moreover, there are $5$ lines
  intersecting $l_1$ and $l_3$; one is $l_2$ and of the $4$ remaining at
  least one, say $l_4$,
  must be fixed. We complete these four lines to an $L$ set and we can
  assume it is $L_b$.
  Again, we study where $L_b$ can be sent by $M$ and we collect the possible values for $M$.
  Similar considerations allow us to obtain other matrices in the cases there are no
  fixed lines for $M$. As in the previous case, we see that the condition that $M(S) =S$
  translates into conditions on the parameters of $S$ which give either
  singular cubic surfaces or smooth cubic surfaces with at least three Eckardt points.\\
  (3) The case in which $\ord(M)=5$ can be solved in a similar way.
\end{proof}

Consider now the Eckardt families $\SE_i$, for $i \in \{1, 2, 3, 4, 6, 9, 10, 18\}$,
introduced in the previous section. 

\begin{definition}
  A permutation $\pi \in \ese$ of the $27$ lines of a surface $S\in \SE_n$
 is called
  \emph{$n$-admissible} if  it maps three coplanar lines passing through an Eckardt point
  to three coplanar lines passing through an Eckardt point.
\end{definition}

 The set of $n$-admissible permutations will be denoted by $\adm_n$ and is a
  subgroup of the Weyl group~$\ese$.
  
Notice that, since the Eckardt points of each surface $S\in \SE_n$ are in a precise
configuration, representable by $n$ formal line triplets,
the above definition does not depend on the particular surface $S$ of $\SE_n$. 

\medskip
Using Proposition~\ref{staine6} and Corollary~\ref{inversapi}, we obtain  the following property.

\begin{prop}
For any $n \in \{1, 2, 3, 4, 6, 9, 10, 18\}$, let $S\in \SE_n$. Then the group monomorphism $\phi$ defined in Proposition~\ref{staine6} restricts to
    \[
    \phi: \st(S) \longrightarrow \adm_n .   \]
The  inverse map (when defined) is given by
    $\pi \mapsto M(L_b, \pi(L_b))$. In particular, $\st(S)$ is contained in the
  set of matrices
  \[
  \mm_n = \{M(L_b, \pi(L_b)) \mid \pi \in \adm_n\}.
  \]
\end{prop}

\begin{remark}
As observed in Section \ref{Sec2}, an $L$-set is not enough to determine
all the labels of the $27$ lines (see Proposition~\ref{2estesi}), but it is necessary to consider the extended $L$-sets
(see Corollary~\ref{corrisp}). This holds for all  $S \in \SE_n$ with $n \ne 9$. \\
Namely,  if $S \in \SE_9$,  the Eckardt point $\ttp_{22} = (E_5, G_2)$ allows us
to uniquely determine the line $E_5$ as the line which intersects $G_2$ in an Eckardt
point. Indeed, $\ttp_{27} = (E_6, G_2)$ is not an Eckardt plane. 
\end{remark}

As a consequence of this observation,
\[
|\adm_n| = 2 |\mm_n|, \; \hbox{for} \; n \ne 9, \quad \hbox{and} \quad |\adm_9| =  |\mm_9| .
\]

The computation of the groups $\adm_n$ is fast, since we have only to
select, among the elements of $\ese$ (permutations of the $27$ simbols $E_1, \dots, F_{56}$ preserving the incidence relations), those that preserve also the Eckardt points.
{}From the elements of $\adm_n$ we can choose the elements of $\mm_n$ (and represent them 
 by simbolic $L$-sets). It turns out that the order of the sets $\mm_n$ are the following:
\[
\begin{array}{llll}
  |\mm_1| = 576 & |\mm_2| = 96 & |\mm_3| = 108 & |\mm_4| = 36 \\
  |\mm_6| = 48 & |\mm_9| = 1296 & |\mm_{10}| = 120 & |\mm_{18}| = 648.
\end{array}
\]
The next step is to explicitly determine  $\st(\SE_n)$: 
for all $L$-sets $L \in \mm_n$, we compute
the matrix $M = M(L_b, L)$. Then $\st(\SE_n)$ is just the collection of all the above
matrices $M$ such that $M(S) = S$ (where $S$ is the generic surface of $\SE_n$).
We get in this way the stabilizers and their orders:
\[
\begin{array}{llll}
  |\st(\SE_1)| = 2 & |\st(\SE_2)| = 4 &|\st(\SE_3)| = 6 &|\st(\SE_4)| = 12 \\
  |\st(\SE_6)| = 24 &|\st(\SE_9)| = 54 &|\st(\SE_{10})| = 120 &
  |\st(\SE_{18})| = 648 
\end{array}
\]

In particular, Lemma~\ref{stabConEck} and the above
computations prove the following
relevant  fact (see also \cite{segre}, Ch.~IV, Sect.~XIV).
\begin{theorem}
  The stabilizer of a smooth cubic surface $S$ is non trivial if and only if
  $S$ has some Eckardt points. 
\end{theorem}

Let us sketch the computation required  in a particular case. 

\begin{example} 
\label{explmat}
The family $\SE_6$ consists of cubic surfaces having six Eckardt points given by
  \[
  (E_1, G_4), (E_2, G_1), (E_2, G_3), (E_3, G_4), (E_4, G_2),
  (F_{13}, F_{24})
  \]
  (see Table~\ref{SEi} and Table~\ref{45planes}). Its equation is given in Table~\ref{listofequ}
  and the corresponding singular locus is  $\Sigma_6 = c(c-e)(3c + e)(c + e)(5c^2 + 2ce + e^2)$. 
  
A direct computation shows that the group $\adm_6$ has $96$ elements, representable  by
   extended $L$-sets, like
     \[
  \begin{array}{c}
   (E_1, G_4, E_2, G_3, E_3, E_5) = \lbe, (E_1, G_4, E_2, G_3, E_3, E_6),
  (E_1, G_4, E_2, F_{12}, E_3, F_{45}),\\
  (E_1, G_4, E_2, F_{12}, E_3, F_{46}), (E_1, G_4, F_{24}, G_2, F_{34}, E_5),
  (E_1, G_4, F_{24}, G_2, F_{34}, E_6),\dots
  \end{array}
  \]
  By deleting the last element of the extended $L$-sets above, we obtain $48$ distinct  $L$-sets which represent the matrices of $\mm_6$.
  For instance, $L = (E_1, G_4, F_{24}, F_{13}, F_{34})$ represents the
  matrix $M = M(L_b, L)  \in \mm_6$, where
  \[
  M=
  \left(
  \begin{array}{cccc}
    c(c+e) & 0 & 0 & (e-c)(3c+e) \\
    0 & c(c+e) & 0 & c^2-4ce-e^2 \\
    0 & 0 & c(c+e) & 2c(c+e) \\
    0 & 0 & 0 & c(c+e)
  \end{array}
\right)
\]
and it is easy to verify that  this matrix belongs to $\st(\SE_6)$. By repeating the same check for all the 48 elements of  $\mm_6$, we see that only $24$  of these matrices stabilize $S \in \SE_6$, so they are the elements of the group $\st(\SE_6)$. 
\end{example}


\section{Subfamilies with larger stabilizer}
\label{Sec5}

So far, we have computed the stabilizers $\st(S)$
for the \emph{generic} cubic surface  $S \in \SE_n$. In this section, 
we check whether, for some $n$, there are  smooth surfaces in $\SE_n$, with $n$ Eckardt points and  a larger stabilizer.

The case of $\SE_1$  requires some attention in order
to avoid hard computations, while subfamilies of the families $\SE_n$ ($n \ge 2$) with larger stabilizer
can be detected in a simpler way.
In these cases, indeed, it is not too difficult to understand if in the group $\adm_n$ there are further $L$-sets that 
 stabilize cubic surfaces under some conditions on the parameters.

  We obtain that
in the families $\SE_i$, for $i \in \{ 2, 4, 6\}$, there are no
cubic surfaces having the same number of Eckardt points and a larger stabilizer.

\subsection {The subfamily  $\SE_1'$.}

If $S\in \SE_1$ contains only one Eckardt point, then, by Lemma~\ref{stabConEck}, the order of $\st(S)$ is  $2^r$, where $1 \le r \le 7$. \\
Following the proof of Lemma~\ref{stabConEck}, we again select the $12$ $L$-sets
which give rise to the $12$ matrices $M_i$. One of them, say $M_1$,  comes from the $L$-set
$(E_1, G_4, F_{24}, F_{13}, F_{35})$. One can see that $\ord(M_1) = 2$ and that $M_1$ stabilizes
$\SE_1$, i.e. $M_1$ generates the cyclic group $\st(S)$, for the general $S$. \\
 At this point we check whether, for some values of the parameters, we can find
surfaces stabilized by other matrices of order 2.
It turns out that the conditions we
get give only singular cubic surfaces, so in $\st(S)$, for all possible $S\in \SE_1$, there is
only one element of order $2$. \\
Therefore, if $\st(S)$ contains other elements, at
least one of them, say $M_2$,  must have order~$4$. In this case, $M_2$ fixes
the plane $\ttp_3 = (E_1, G_4, F_{14})$ and also, consequentely, the three lines of $\ttp_3$.
Since $E_2$ meets $G_4$, then $M_2(E_2)$ is  another line which intersects $G_4$; it is easy to see that it can be assumed to be $E_3$. \\
 Hence we collect all the $L$-sets of the form $(E_1, G_4, E_3, *, *)$ and
check if they give matrices which, under suitable conditions on the parameters,
stabilize the corresponding cubic surface. In this way, we find that the matrix given by
the $L$-set $(E_1, G_4, E_3, G_6, F_{24})$ is in the stabilizer
of a cubic surface
of $\SE_1$ as long as the conditions:
\begin{equation}
\begin{array}{r@{}l}
  e &{} =  \displaystyle \frac{(2\sqrt{-1} - 1)c(5c - (4\sqrt{-1} - 3)f)}
  {5(c-f)} \\
  d &{} =  \displaystyle
  \frac{(2\sqrt{-1} + 1)(5c^2 - (4\sqrt{-1} + 8)cf - 5f^2)}{5(f-c)}  
\end{array}
\label{condSe1'}
\end{equation}
are satisfied (note that the denominator is not zero, since the polynomial
$c-f$ is a factor of $\Sigma_1$). \\
The conditions~(\ref{condSe1'}) plus condition $\Lambda_1$ of
Table~\ref{SEi} will be denoted by $\Lambda_1'$ and  give a
family $\SE_1'$ depending on two parameters ($c$ and $f$) such that
$\st(\SE_1')$ is (in general) of order $4$ (and has only one Eckardt point).
Other $L$-sets give other families of cubic surfaces with a matrix of order
$4$ in the stabilizer, but are all projectively equivalent to  $\SE_1'$.

\subsection {The subfamily  $\SE_1''$.}

Now  we impose to the generic cubic of $\SE_1'$ that the
matrix corresponding to the $L$-set $(G_4, E_1, F_{15}, E_5, F_{12})$
is in the stabilizer, we get the conditions
\begin{equation}
\begin{array}{r@{}l}
  c &{} = \sqrt{-1} - \sqrt{2}\\
  f &{} = 3
\end{array}
\label{condSe1''}
\end{equation}
which gives a cubic surface $\SE_1''$
with one Eckardt point and stabilizer of order $8$.
Also the matrices obtained from other $L$-sets, like 
\[
(E_1, G_4, E_2, F_{12}, E_3),\quad (E_1, G_4, F_{46}, G_6, F_{45}),\quad
(E_1, G_4, E_6, G_5, E_5)
\]
and several others,
give cubic surfaces with
stabilizer of order $8$, but these surfaces are all projectively equivalent
to $\SE_1''$. Therefore the condition 
$\Lambda_1''$ given by $\Lambda_1'$ plus conditions~\ref{condSe1''}
(when substitute
into the polynomial (\ref{cubicaB})), gives the unique
surface with one Eckardt point and stabilizer of order $8$.

\begin{remark}
  In the list of possible stabilizers of cubic surfaces given in~\cite{segre}
  this case is missing, as also remarked in~\cite{hosoh}.
\end{remark}

\subsection {The subfamily  $\SE_9'$.}

Finally, the family $\SE_9$ contains cubic surfaces with larger stabilizer. One
of them is  obtained from $\SE_9$ by imposing the conditions:
\begin{equation}
\begin{array}{r@{}l}
  f&{} =  1/4(\sqrt{-1}-\sqrt{3})^3 -(\sqrt{-1}-\sqrt{3})-1\\
  c&{} = 1
\end{array}
\label{condSe9'}
\end{equation}
on its parameters. Hence the conditions $\Lambda_9'$, obtained by adding (\ref{condSe9'}) to
$\Lambda_9$,  give the unique (up to projectivities)
cubic surface, say $\SE_9'$, still having $9$ Eckardt points, but stabilizer of order $108$.


\section{Structure of the automorphisms groups}
\label{Sec6}

The knowledge of the explicit equations of the cubic surfaces of each of
the Eckardt families $\SE_i$, their lines and their stabilizers allows us to 
determine the structure of these groups and also a  graphic representation of them.

Here we first describe the stabilizers
giving (in a rough way) for each group, the generators in terms of permutations of
the numbers $1,\dots, 27$ which represent the $27$ lines, according to the
following correspondence:

{\small
\[
\begin{array}{cccccccccccccccc}
  E_1& E_2& E_3& E_4& E_5& E_6 & G_1& G_2& G_3& G_4& G_5& G_6 &&&&\\ \bigskip
  1 & 2 & 3 & 4 & 5 & 6 &7 & 8 & 9 & 10 & 11 & 12&&&&\\ 
  F_{12}& F_{13}& F_{14}& F_{15}& F_{16}& F_{23}& F_{24}& F_{25}& F_{26}& F_{34}& F_{35}& F_{36}& F_{45}& F_{46}& F_{56}\\
  13 & 14 & 15 & 16 & 17 & 18 & 19 & 20 & 21 & 22 & 23 & 24 & 25 & 26 & 27 
\end{array}
\]
}


\begin{table}[H]
  \caption{Stabilizers of Eckardt families and their generators in terms of permutations of the lines.}
  \label{stabitutti}
  \small{
  \begin{tabular}{|c|l|l|c|} \hline
    $S$ & type of $\st(S)$ & generators & $|\st(S)|$ \\ \hline
    $\SE_1$ & $C_2$ & $g_1$ &  2 \\
     $\SE_1'$ & $C_4$  & $g_1'$ & 4 \\
     $\SE_1''$ & $C_8$  & $g_1''$ & 8 \\
\hline
      $\SE_2$ & $C_2 \times C_2$  & $g_2, h_2$ & 4 \\
    $\SE_3$ & $S_3$  & $g_3, h_3$ & 6 \\
    $\SE_4$ & $C_2 \times S_3$  & $g_4, h_4$ & 12 \\
\hline
    $\SE_6$ & $S_4$  & $g_6, h_6$ & 24 \\
\hline
    $\SE_9$ & $((C_3 \times C_3) \rtimes C_3)\rtimes C_2$  & $g_9, h_9, k_9$ & 54 \\     
    $\SE_9'$ & $((C_3 \times C_3) \rtimes C_3)\rtimes C_4$  & $g_9', h_9'$ & 108 \\
\hline
    $\SE_{10}$ & $S_5$  & $g_{10}, h_{10}$ & 120  \\
\hline
    $\SE_{18}$ & $(C_3 \times C_3 \times C_3)\rtimes S_4$  &$g_{18}, h_{18}$ & 648  \\
\hline
\end{tabular}
  }
\end{table}

where:
{\scriptsize
\begin{eqnarray*}
g_1 & = &  (2,19)(3,22)(4,7)(5,25)(6,26)(8,13)(9,14)(11,16)(12,17)(18,27)(20,24)(21,23)\\
g_1' & = &  (2,3,19,22)(4,21,7,23)(5,26,25,6)(8,16,13,11)(9,12,14,17)(18,24,27,20)\\
g_1'' & = &  (1,10)(2,11,3,8,19,16,22,13)(4,24,21,27,7,20,23,18)(5,17,26,9,25,12,6,14)\\
g_2 & = &  (1,15)(3,22)(4,8)(5,25)(6,26)(7,13)(9,18)(11,20)(12,21)(14,27)(16,24)(17,23)\\
h_2 & = &  (2,19)(3,22)(4,7)(5,25)(6,26)(8,13)(9,14)(11,16)(12,17)(18,27)(20,24)(21,23)\\
g_3 & = &  (2,19)(3,22)(4,7)(5,25)(6,26)(8,13)(9,14)(11,16)(12,17)(18,27)(20,24)(21,23)\\
h_3 & = &  (1,9,14)(2,19,10)(3,13,4)(5,21,17)(6,20,16)(7,8,22)(11,24,26)(12,23,25)(15,18,27)\\
g_4 & = &  (2,19)(3,22)(4,7)(5,25)(6,26)(8,13)(9,14)(11,16)(12,17)(18,27)(20,24)(21,23)\\
h_4 & = &  (1,18,14,15,9,27)(2,19,10)(3,7,4,22,13,8)(5,12,17,25,21,23)(6,11,16,26,20,24)\\
g_6 & = &  (1,13)(2,10)(3,18)(5,20)(6,21)(7,15)(9,22)(11,25)(12,26)(14,27)(16,24)(17,23)\\
h_6 & = &  (1,3,15,22)(2,19)(4,13,14,18)(5,6)(7,27,9,8)(11,23,21,16)(12,24,20,17)(25,26)\\
g_9 & = &  (2,19)(3,22)(4,7)(5,25)(6,26)(8,13)(9,14)(11,16)(12,17)(18,27)(20,24)(21,23)\\
h_9 & = &  (1,9,14)(2,19,10)(3,13,4)(5,21,17)(6,20,16)(7,8,22)(11,24,26)(12,23,25)(15,18,27)\\
k_9 & = &  (1,21,3)(2,7,13)(4,27,20)(5,26,19)(6,24,9)(8,12,14)(10,16,23)(11,15,22)(17,25,18)\\
g_9' & = &  (2,3,19,22)(4,20,7,24)(5,6,25,26)(8,17,13,12)(9,11,14,16)(18,23,27,21)\\
h_9' & = &  (1,4,6)(2,3,5)(7,10,12)(8,9,11)(13,22,27)(14,25,21)(15,26,17)(16,19,24)(18,23,20)\\
g_{10} & = &  (1,13)(2,9)(4,19)(5,20)(6,21)(7,14)(10,22)(11,23)(12,24)(15,27)(16,26)(17,25)\\
h_{10} & = &  (1,3,15,8,7)(2,9,22,27,19)(4,13,14,10,18)(5,17,11,24,21)(6,16,12,23,20)\\
g_{18} & = &  (1,17)(2,5)(3,4,24,26)(6,14,22,15)(7,9,20,10)(8,19,16,18)(11,25,13,23)(21,27)\\
h_{18} & = &  (1,26,7)(2,11,20)(3,8,18)(4,17,10)(5,9,23)(6,15,12)(13,16,14)(19,25,22)(21,27,24)
\end{eqnarray*}
}

\begin{remark} The content of the above Table clearly coincides with  Table 9.6 in  ~\cite{dol} and  with the description of the groups given in 
~\cite{hosoh}.
  Our approach  is, however,
  totally explicit and allows us to immediately obtain the projectivities (i.e.\ the matrices)
  which stabilize the surfaces, accordingly to the procedure in Example~\ref{explmat}.
\end{remark}

There are however other ways to interpret the stabilizers which show hidden  symmetries involving Eckardt points, tritangent planes, Sylvester pentahedron, \dots

\subsection{Stabilizers of $\SE_1$, $\SE_1'$ and $\SE_1''$}

\begin{itemize}

\item  If $S \in \SE_1$ is generic, then   $\st(S) \cong C_2$ acts as follows:

- the three lines $E_1$, $G_4$, $F_{14}$ which define the Eckardt point are fixed;

-  the five tritangent planes through $E_1$ (respectively, through $G_4$  and through $F_{14}$) are fixed;

- the non trivial element of $\st(S)$ exchanges the two lines
$E_2$ and $F_{24}$, so $\st(S)$ is isomorphic to the subgroup of $\ese$ generated by the
$2$-cycle $(E_2, F_{24})$ (or, analogously, by the $2$-cycle
$(\ttp_{12}, \ttp_{31})$ if we consider the action of $\st(S)$ on the planes).

 \item If $S \in \SE_1'$ is generic, then  $\st(S) \cong C_4$ acts as follows:
 
- the three lines $E_1$, $G_4$, $F_{14}$ are fixed;

- the only tritangent fixed plane is the Eckardt plane $\ttp_3$.

-  the group $\st(S)$ is isomorphic to
the cyclic group generated by the $4$-cycle $(E_2, F_{34}, F_{24}, E_3)$
(or by the $4$-cycle of planes $(\ttp_{9}, \ttp_{31}, \ttp_{41}, \ttp_{12})$).

\item If $S = \SE_1''$, then  $\st(S) \cong C_8$. The only  fixed line is $F_{14}$ and that 
$\st(S)$  is isomorphic to
the group generated by the $8$-cycle
\[
(E_2, G_2, F_{34}, G_5, F_{24}, F_{12}, E_3, F_{15})
\]
(or by the $8$-cycle of planes
$(\ttp_9, \ttp_{17}, \ttp_{31}, \ttp_{14}, \ttp_{41},
\ttp_{6}, \ttp_{12}, \ttp_{42})$).

\end{itemize}

\begin{figure}
  \includegraphics[height= 4cm]{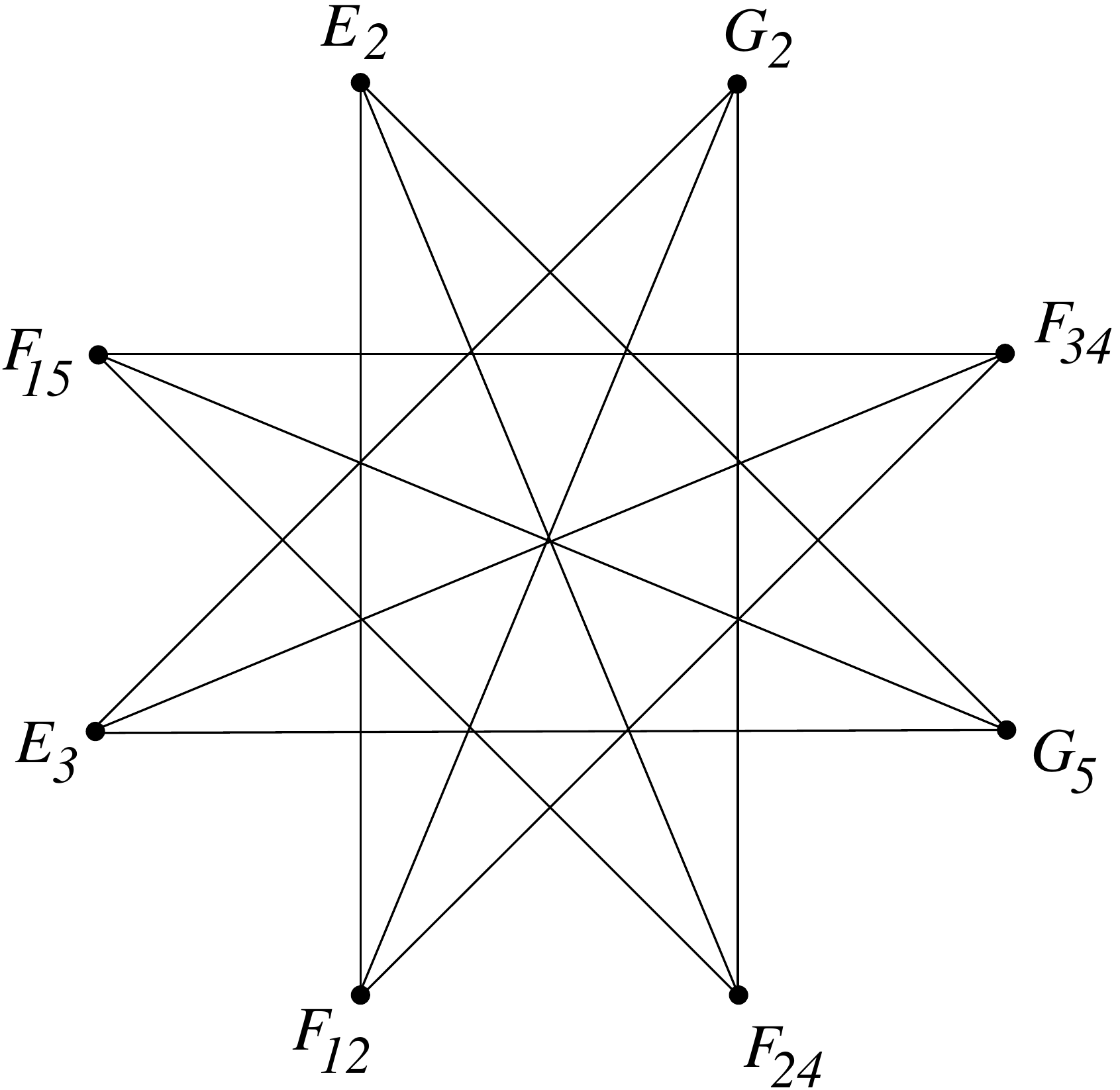}
  \caption{Representation of the stabilizers of $\SE_1, \SE_1'$ and $\SE_1''$.
    Two vertices are connected by a line if the corresponding lines are
    incident. The groups $\st(S)$ for $S$ in $\SE_1, \SE_1', \SE_1''$
    are the rotations around the center of the figure by $180^o$, by $90^o$
    and by $45^o$, respectively.
  \label{octagon}}
\end{figure}

\medskip
\noindent
Figure~\ref{octagon} allows us to visualize these three stabilizers. 
The  $8$ vertices $E_2$, $G_2$, $F_{34}, \dots$ 
represent $8$ lines of $S$ and two vertices are connected if and only
if the corresponding lines are coplanar (and hence determine a tritangent
plane). For instance, the vertices $E_2$ and $G_5$ are connected,
since the corresponding lines meet on the tritangent plane $\ttp_9 = (E_2, G_5, F_{25})$.

In the stabilizer of a generic element  $S \in \SE_1$ we have the $2$-cycles $(E_2, F_{24})$, 
$(G_2, F_{12})$, $(F_{34}, E_3)$ and  $(G_5, F_{15})$, which correspond to the rotation of $180^o$ of Figure~\ref{octagon}  around its center. Hence 
$\st(S)$ can be seen as a rotation group.

The action of $\st(S)$ on the tritangent planes
fixes the planes $\ttp_8 = (E_2, F_{24})$, $\ttp_1$, $\ttp_{13}$,
$\ttp_4$ and exchanges the planes of the couples $(\ttp_6, \ttp_{17})$,
$(\ttp_{12}, \ttp_{31})$, $(\ttp_{42}, \ttp_{14})$,
$(\ttp_9, \ttp_{41})$ and these permutations are also coherent with
the rotation of the Figure~\ref{octagon}.

Similarly, the stabilizer of a generic cubic of $\SE_1'$ is
isomorphic to the group 
generated by the clockwise rotation of the figure around its center
by $90^o$. Also here, a rotation of the figure can be identified
either with a permutation of lines or of tritangent planes. 

Finally, the stabilizer of $\SE_1''$ is isomorphic to the group
generated by
the clockwise rotation
of the figure around its center of $45^o$.

\subsection{Stabilizer of $\SE_2$}
The generic smooth surface $S \in \SE_2$ has  two Eckardt points $(E_1,G_4), (E_2,G_4)$
(see Remark~\ref{osservase2}). The only  line fixed by $\st(S)$ is $G_4$.
The five tritangent planes through $G_4$ are fixed as well and, so, also
the two Eckardt points.

If  $\{e, g_1, g_2,g_3\}$ are  the four elements of the stabilizer and $\{E_1, F_{14}, E_2, F_{24}\}$ are all the lines on the two Eckardt planes, out of $G_4$,  then, up to renaming the elements,  $g_1$ is the  $2$-cycle $(E_1, F_{14})$, $g_2$ is the  $2$-cycle $(E_2, F_{24})$ and 
 $g_3$ exchanges $E_1$ with $F_{14}$ and $E_2$ with $F_{24}$, i.e. $g_3 = g_1 \circ g_2$.
 This is another description of $\st(S)$  as  $C_2\times C_2$. Moreover, concerning the action of $\st(S)$  on the tritangent planes, it acts
(for instance) on the set of four planes $\{\ttp_1, \ttp_{16}, \ttp_{6}, \ttp_{17} \}$ exchanging
the first two or the last two or both couples.

\subsection{Stabilizer of $\SE_3$}
For a generic $S \in \SE_3$,  the group $\st(S) \cong S_3$  does not fix lines and permutes, in all the possible
ways, the three lines $E_1, G_3, F_{13}$ on the plane
$\ttp_2$. In addition, also the lines on the planes $\ttp_8$ and $\ttp_{37}$ are permuted in all the possible ways. 

The action of $\st(S)$ on the tritangent planes is even more explicit, since corresponds to all the permutations
of the planes $\ttp_3, \ttp_7, \ttp_{34}$, i.e.\ to all the permutations of the
three Eckardt points.

\subsection{Stabilizer of $\SE_4$}
The generic cubic surface $S \in \SE_4$ has  four Eckardt points: $\ttp_8 = (E_2,G_4,F_{24})$, 
$\ttp_3 = (E_1, G_4, F_{14})$,
$\ttp_7 = (E_2, G_3, F_{23})$ and $\ttp_{34} = (F_{13}, F_{24}, F_{56})$.
They belong to the Eckardt plane $\ttp_8$ and three of
them are collinear: $\ttp_3,\ttp_7,\ttp_{34} $. Clearly, each of these three points belongs to one of the three lines through $\ttp_8$
and this plane is fixed by $G = \st(S)$.

The action of $G$ on the lines gives rise to the following five orbits:
\[
\begin{array}{l}
\{E_1, F_{23}, F_{13}, F_{14}, G_3, F_{56}\}, \quad
\{E_2, F_{24}, G_4\}, \quad
\{E_3, F_{34}, G_1, F_{12}, G_2, E_4\}, \\
\{E_5, F_{45}, G_6, F_{26}, F_{35}, F_{16}\}, \quad
\{E_6, F_{46}, G_5, F_{25}, F_{36}, F_{15}\}
\end{array}
\]
The action of $G$ on the three
collinear Eckardt points $\ttp_3$, $\ttp_7$, $\ttp_{34}$ gives all their
permutations, meanwhile $G$ exchanges the planes $\ttp_2$
and $\ttp_{37}$. This is another way to see that  $G$ is isomorphic to the direct product
$C_2 \times S_3$.

The group $G$ can be also represented by considering an hexagon whose vertices are the six lines of  one of  the orbits above. Rearranging the labels of the vertices in a suitable way, $G$  turns out to be isomorphic to the dihedral group $D_6$.

\subsection{Stabilizer of $\SE_6$}
\begin{figure}
  \includegraphics[height= 4cm]{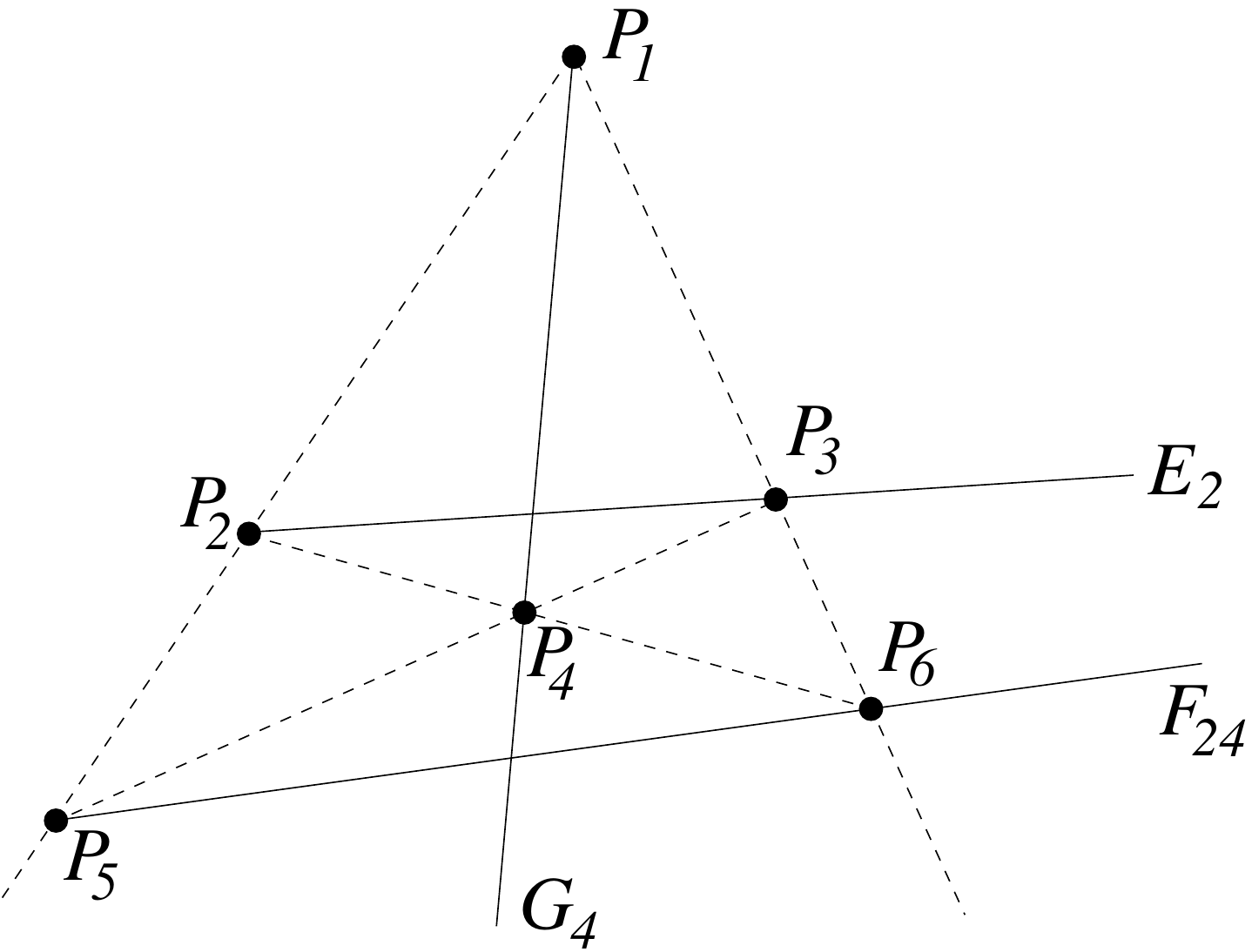}
  \caption{Reciprocal position of the six Eckardt points of a cubic of
    $\SE_6$
  \label{figA}}
\end{figure}

The group $G = \st(S)$ for a generic $S \in \SE_6$ is $S_4$.
The six Eckardt points  are $P_1= \ttp_3, P_2 = \ttp_6, P_3 = \ttp_7,
P_4 = \ttp_{13},
P_5 = \ttp_{17}, P_6 = \ttp_{34}$ and are all contained in the plane $\ttp_8$ of equation $x=0$. \\
The action of $G$ on the tritangent planes has seven orbits, two of them
are:
\[
\{\ttp_1, \ttp_{11}, \ttp_{18}, \ttp_{37} \},
\quad \{\ttp_2, \ttp_{12}, \ttp_{16}, \ttp_{31} \}.
\]
The action of $G$, restricted on each of these two orbits, gives all the
possible permutations of  the four planes, confirming that $G \cong S_4$. \\
Moreover it is straightforward to see that
\[
\begin{array}{ll}
  S \cap \ttp_2 \cap \ttp_{37} = \{P_1, P_3, P_6  \}, &
  S \cap \ttp_{11} \cap \ttp_{31} = \{P_2, P_4, P_6  \}, \\
  S \cap \ttp_{12} \cap \ttp_{18} = \{P_3, P_4, P_5  \}, &
  S \cap \ttp_{1} \cap \ttp_{16} = \{P_1, P_2, P_5  \}.
\end{array}
\]
Hence  $G$ can also be seen as the group of permutations of the dotted lines in Figure~\ref{figA}: they  are not lines of $S$ 
and each of them contains three Eckardt points. 

\subsection{Stabilizer of $\SE_9$}
For a generic $S \in \SE_9$, the nine Eckardt points ($\ttp_3$, $\ttp_7$, $\ttp_{14}$, $\ttp_{20}$,
$\ttp_{22}$,
$\ttp_{26}$, $\ttp_{33}$, $\ttp_{34}$, $\ttp_{42}$)
are coplanar, belonging to the plane $\pi: (1-\sqrt {-3}) x -y +z=0$. Moreover,   they are the nine inflection points of
the cubic curve $\pi \cap S$ and are therefore in the Hesse configuration.

\begin{figure}
  \includegraphics[height= 6.5cm]{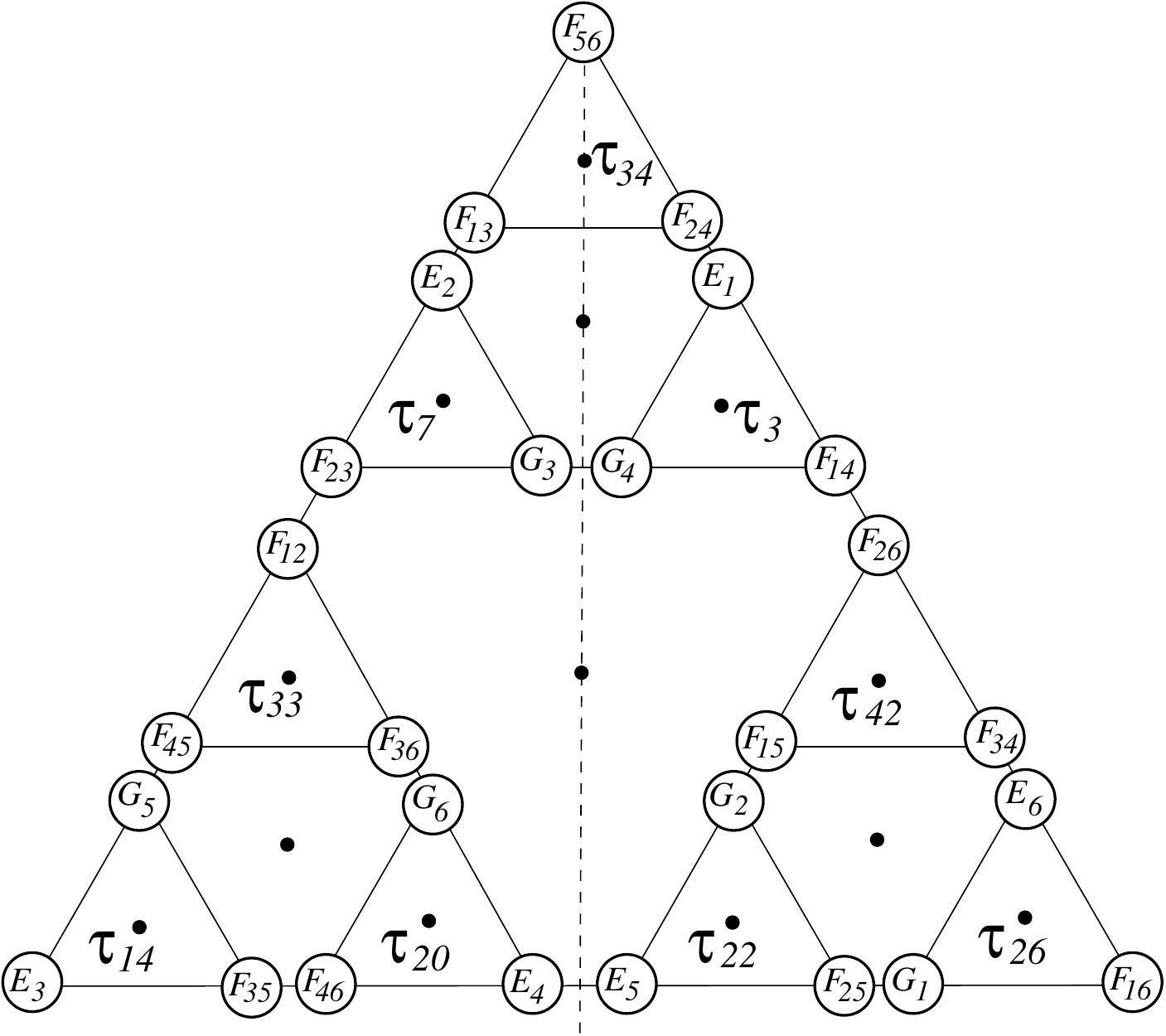}
  \caption{Group of symmetries of $\SE_9$ and
    of $\SE_9'$. Each small triangle represents an Eckardt point,
    for instance the
    triangle whose vertices are $F_{13}, F_{24}, F_{56}$ is the Eckardt pont
    $\ttp_{34}$ (see  Table~\ref{45planes}).
  \label{fig9Pt}}
\end{figure}

The group $\st(S)$ contains the following permutations:
\begin{itemize}
\item  $g_1$,  given by the product of the disjoint cycles
  \[
  \begin{array}{l}
  (E_{1}, F_{14}, G_{4}), (E_{2}, G_{3}, F_{23}), (E_{3}, G_{5}, F_{35}),
    (E_{4}, F_{46}, G_{6}), (E_{5}, G_{2}, F_{25}), \\
    (E_{6}, F_{16}, G_{1}),
(F_{12}, F_{36}, F_{45}), (F_{13}, F_{56}, F_{24}), (F_{15}, F_{26}, F_{34}),
  \end{array}
  \]
 fixes the $9$ Eckardt points and rotates clockwise around their
  baricenters the small triangles of Figure~\ref{fig9Pt} of $120^o$;
\item $g_2$,  given by the product of the disjoint cycles
  \[
  \begin{array}{c}
  (E_{1}, G_{3}, F_{13}), (E_{2}, F_{24}, G_{4}), (E_{3}, F_{12}, E_{4}),
 (E_{5}, F_{26}, F_{16}), (E_{6}, F_{25}, F_{15}), \\(G_{1}, G_{2}, F_{34}),
 (G_{5}, F_{36}, F_{46}), (G_{6}, F_{35}, F_{45}),(F_{14}, F_{23}, F_{56}),
  \end{array}
  \]
 rotates clockwise around their baricenters the three
  medium triangles of Figure~\ref{fig9Pt} of $120^o$;
\item $g_3$,  given by the product of the disjoint cycles
  \[
  \begin{array}{c}
 (E_{1}, G_{2}, F_{12}), (E_{2}, F_{15}, E_{4}), (E_{3}, F_{56}, F_{16}),
    (E_{5}, F_{45}, G_{4}), (E_{6}, F_{35}, F_{13}), \\
    (G_{1}, G_{5}, F_{24}),
 (G_{3}, F_{26}, F_{46}), (G_{6}, F_{23}, F_{34}),(F_{14}, F_{25}, F_{36}),
  \end{array}
  \]
rotates clockwise around its baricenter the large triangle
  of Figure~\ref{fig9Pt} and, simultaneously, rotates clockwise the three
  small triangles $(E_1, F_{14}, G_4)$, $(G_2, F_{25}, E_5)$ and
  $(F_{12}, F_{36}, F_{45})$ of $240^o$ and the tree small triangles
  $(E_2, G_3, F_{23})$, $(F_{15}, F_{26}, F_{34})$, $G_6, E_4, F_{46})$ of
  $120^o$; therefore also $g_3$ has order $3$.
\item $g_4$, given by the product of disjoint cycles
  \[
  \begin{array}{c}
  (E_{1}, G_{3}), (E_{2}, G_{4}), (E_{3}, G_{1}), (E_{4}, G_{2}),
(E_{5}, G_{6}), (E_{6}, G_{5}),\\
(F_{12}, F_{34}), (F_{14}, F_{23}), (F_{15}, F_{36}), (F_{16}, F_{35}),
(F_{25}, F_{46}), (F_{26}, F_{45})
  \end{array}
  \]
 is the reflection of the large triangle of Figure~\ref{fig9Pt}
  along the dotted line.
\end{itemize}
This shows that the group $\st(S)$ is isomorphic to
\[
\left((\langle g_1\rangle \times \langle g_2\rangle) \rtimes \langle g_3
\rangle\right)\rtimes \langle g_4 \rangle \cong ((C_3 \times C_3) \rtimes C_3)\rtimes C_2 .
\]

\subsection{Stabilizer of $\SE_9'$}
Since $\SE_9' \in \SE_9$, then $\st(\SE_9)$ is a subgroup of $\st(\SE_9')$ (see Sections~\ref{Sec4} and  \ref{Sec5}).

Obviously, the nine Eckardt points are again coplanar and are the nine inflection
points of the cubic curve obtained as intersection of their plane and
the surface. 

Nevertheless, the group structure of $\st(\SE_9')$ is different. 
Keeping the notation for  $g_1, g_2, g_3, g_4$, defined in the case $\SE_9$, let us set
 $g_5$ to be  the permutation of the lines
given by the product of the following disjoint cycles:
\[
\begin{array}{l}
  (E_1, F_{25}, G_3, F_{46}), (E_2, E_{4}, G_4, G_2), (E_3, F_{26}, G_1, F_{45}),\\
(E_5, F_{23}, G_6, F_{14}), (E_6, F_{12}, G_5, F_{34}),
  (F_{15}, F_{16}, F_{36}, F_{35}).
  \end{array}
\]
It is not difficult to see that $\ord(g_5) = 4$ and  $g_5^2 = g_4$. Moreover, $g_5$ acts
on the small triangles of Figure~\ref{fig9Pt} (i.e.\ on the nine
Eckardt points) as follows:
\[
\ttp_3 \mapsto \ttp_{22} \mapsto \ttp_7 \mapsto \ttp_{20} \mapsto \ttp_3, \quad
\ttp_{14} \mapsto \ttp_{42}\mapsto \ttp_{26} \mapsto \ttp_{33} \mapsto \ttp_{14},
\quad \ttp_{34} \mapsto \ttp_{34}
\]
Therefore the group structure of the stabilizer is
\[
\st(\SE_9') = 
\left((\langle g_1\rangle \times \langle g_2\rangle) \rtimes
\langle g_3\rangle\right) \rtimes \langle g_5\rangle \cong ((C_3 \times C_3) \rtimes C_3)\rtimes C_4.
\]

\subsection{Stabilizer of $\SE_{10}$} 
\begin{figure}
  \includegraphics[height= 5cm]{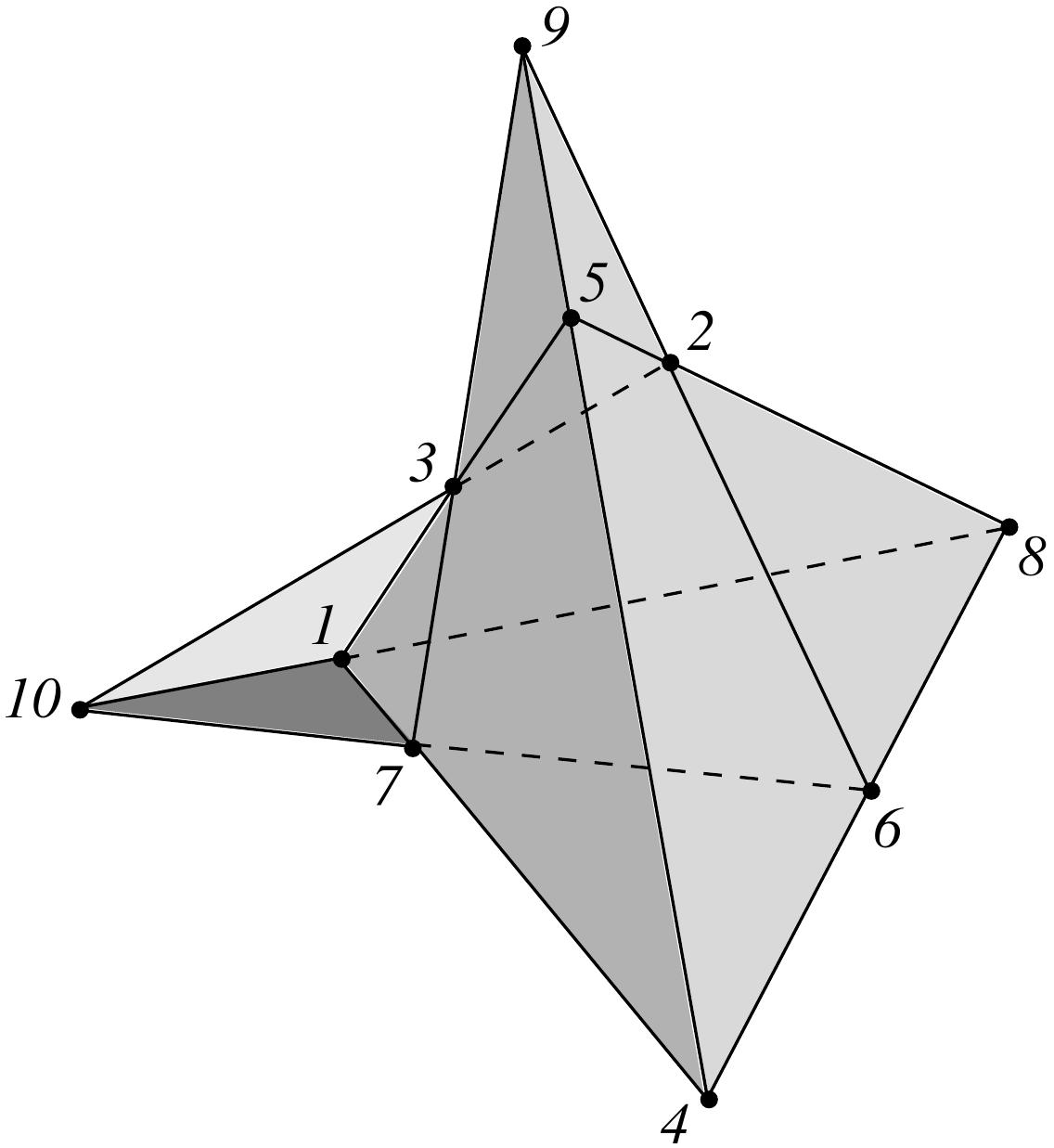}
  \caption{The Sylvester pentahedron of $\SE_{10}$.
    The five planes $\ttp_2$, $\ttp_6$,
    $\ttp_{13}$, $\ttp_{17}$ and $\ttp_{37}$  are the planes which
    contain, respectively,
    the Eckard points: $(1, 2, 3, 5, 8, 10)$, $(1, 3, 4, 5, 7, 9)$,
    $(2, 4, 5, 6, 8, 9)$, $(1, 4, 6, 7, 8, 10)$ and $(2, 3, 6, 7, 9, 10)$.}
  \label{sylv}
\end{figure}

The cubic surface $\SE_{10}$ has the $10$ Eckardt points: $\ttp_1$, $\ttp_3$, $\ttp_7$, $\ttp_8$, $\ttp_{11}$, $\ttp_{12}$,
$\ttp_{16}$, $\ttp_{18}$, $\ttp_{31}$, $\ttp_{34}$. \\
It turns out that the action of $\st(\SE_{10})$ on the $45$ tritangent planes gives rise to an orbit consisting of  five planes, 
$\ttp_2, \ttp_6, \ttp_{13}, \ttp_{17}, \ttp_{37}$  (none of them is an Eckardt plane) and that $\st(\SE_{10})$ permutes them in all possible ways. In this way, we obtain an explicit representation of $\st(\SE_{10})$ as $S_5$. \\
The polynomials defining the above planes  are:
\[
 \ttp_2 : \; y-z,  \quad \ttp_6 : \; x-y, \quad \ttp_{13}  : \;  x+t , \quad  \ttp_{17} : \; x-2y+t , \quad \ttp_{37}  : \;  2x-y+z
\]
and it is easy to verify that $\SE_{10}$ is  the sum of the cubes of these five linear polynomials. Therefore the five planes above compose the Sylvester pentahedron of  
$\SE_{10}$ (see~\cite{segre}, Ch. IV, §84). 
Figure~\ref{sylv} shows the five
planes and the $10$ Eckardt points (there labelled by $1, 2, \dots, 10$).

\subsection{Stabilizer of $\SE_{18}$} 
\begin{figure}[H]
  \includegraphics[height= 5.5cm]{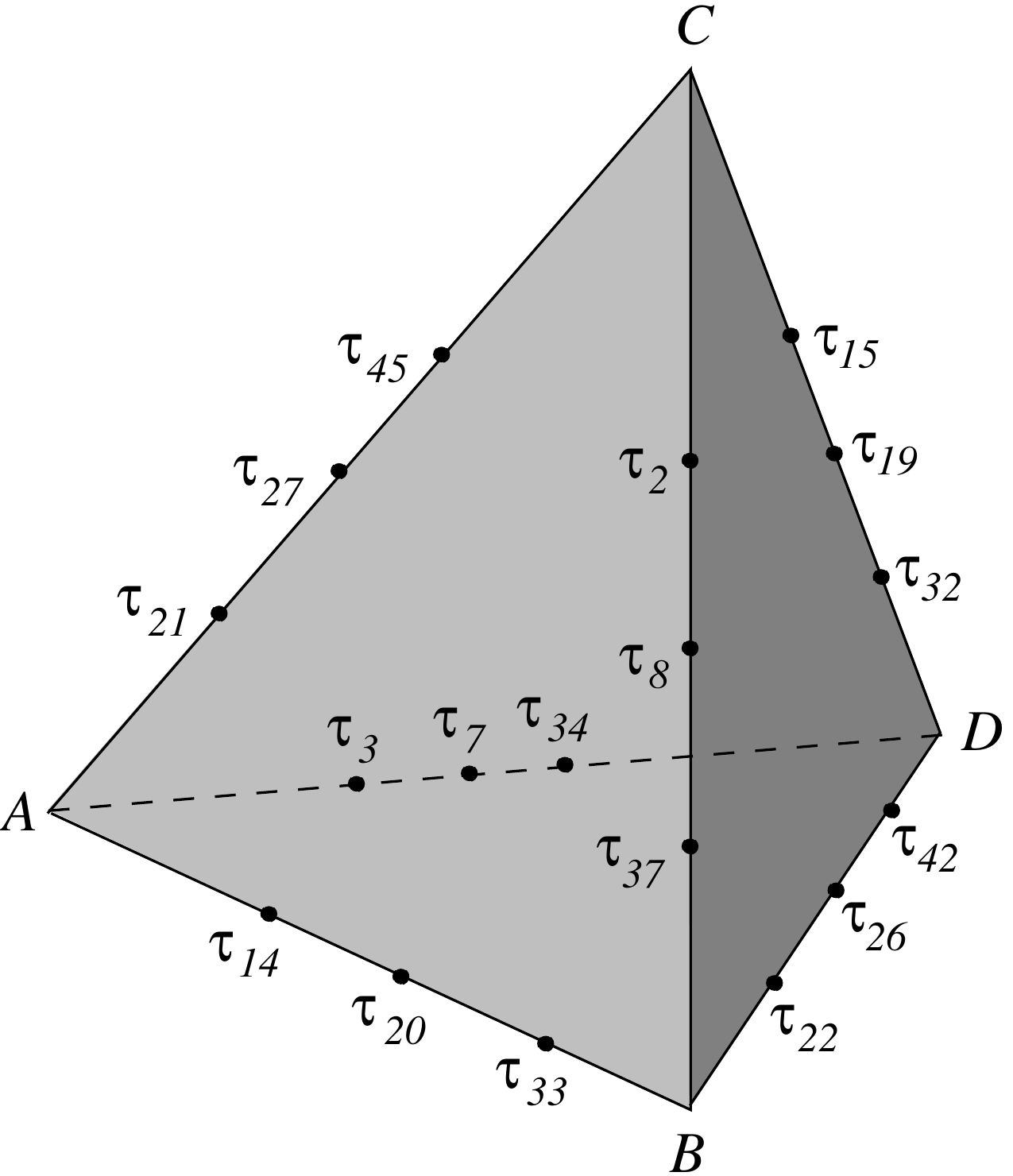}
  \caption{A representation of the group $\st(\SE_{18})$.}
  \label{simm648}
\end{figure}

The cubic surface $\SE_{18}$ has the $18$ Eckardt points: 
\[
\ttp_{2}, \ttp_{3}, \ttp_{7}, \ttp_{8}, \ttp_{14}, \ttp_{15},
\ttp_{19}, \ttp_{20}, \ttp_{21}, 
\ttp_{22}, \ttp_{26}, \ttp_{27},
\ttp_{32}, \ttp_{33}, \ttp_{34}, \ttp_{37}, \ttp_{42}, \ttp_{45}
\]
and  the group $G = \st(\SE_{18})$ has order~$648$. \\
Also in this case, 
the action of $G$ on the $45$ tritangent planes has a particular orbit but, here, it consists of the 
 $18$ Eckardt planes listed above. \\
 A more detailed analysis shows that the $18$ Eckardt points are contained in four planes, represented in 
 Figure~\ref{simm648} as the faces  of a tetrahedron:
 \[
 \pi_1 = A+B+C, \;
\pi_2 = B+C+D, \;
\pi_3 = A+B+D, \;
\pi_4 = A + C + D
\]
and whose equations are:
\[
\begin{array}{ll}
  \pi_1 :   x-(1+\sqrt{-3})y -t, & \quad \pi_2 : (1+\sqrt{-3})(x-t) -4y,\\
  \pi_3 : (1-\sqrt{-3})x-y+z, & \quad \pi_4 : 2(\sqrt{-3}-1)x-(1+\sqrt{-3})(y-z).
\end{array}
\]
Each plane $\pi_i$ contains $9$ Eckardt points, as shown in Figure~\ref{simm648}, and these are precisely 
 the inflection points of the cubic curve $\pi_i \cap \SE_{18}$. \\
The action of $G$ on the $18$ Eckardt points gives a group $H$ of
permutations of these
points and  $|H|= 648$, so  $H \cong G$.
It contains the three permutations:
\[
\begin{array}{l}
  g_1 = (\ttp_3, \ttp_{34},\ttp_{7})(\ttp_{15},\ttp_{32},\ttp_{19})
  (\ttp_{22},\ttp_{42},\ttp_{26})\\
  g_2 = (\ttp_3,\ttp_{34},\ttp_{7})(\ttp_{14},\ttp_{33},\ttp_{20})
  (\ttp_{21},\ttp_{27},\ttp_{45})\\
  g_3 = (\ttp_2,\ttp_8,\ttp_{37})(\ttp_{15},\ttp_{32},\ttp_{19})
  (\ttp_{21},\ttp_{45},\ttp_{27})
  \end{array}.
\]
Each of them acts cyclically on the triplets of collinear Eckardt points, as can be seen in Figure~\ref{simm648}; therefore, setting $N = \langle g_1, g_2, g_3 \rangle$, it is clear that 
$N \cong C_3 \times C_3 \times C_3$. Moreover $N$ is a normal subgroup of $H$ and 
the quotient group $H/N$ is generated by the two classes $[g_4]$ and $[g_5]$,
where:
\[
\begin{array}{l}
  g_4 = (\ttp_2,\ttp_{21})(\ttp_3,\ttp_{42})(\ttp_{7},\ttp_{22})
  (\ttp_{8},\ttp_{45})(\ttp_{20},\ttp_{33})(\ttp_{26},\ttp_{34})
  (\ttp_{27},\ttp_{37})\\
  g_5 = (\ttp_2,\ttp_7)(\ttp_3,\ttp_{37},\ttp_{34},\ttp_{8})
(\ttp_{14},\ttp_{22},\ttp_{32},\ttp_{21})
(\ttp_{15},\ttp_{45},\ttp_{33},\ttp_{42})
(\ttp_{19},\ttp_{27},\ttp_{20},\ttp_{26})
\end{array}.
\]
Concerning the action of these two elements on the four planes of the tetrahedron, note that
\[
g_4: \quad \pi_2 \leftrightarrow \pi_4, \qquad
g_5: \quad \pi_1 \mapsto \pi_3 \mapsto \pi_2 \mapsto \pi_4 \mapsto \pi_1
\]
Therefore $\ord(g_4) = 2$, $\ord(g_5) = 4$ and  $H/N \cong S_4$. Setting  $K = \langle g_4, g_5 \rangle$
the corresponding subgroup  of $H$, we get that $H = N \rtimes K$, hence $\st(\SE_{18})$ is now described by
\[
(\langle g_1 \rangle \times \langle g_2 \rangle \times \langle g_3 \rangle)
\rtimes \langle g_4, g_5 \rangle \cong
(C_3 \times C_3 \times C_3) \rtimes S_4.
\]

\begin{center}
\textbf{Acknowledgement}
\end{center}

We would like to thank Igor Dolgachev for many interesting remarks and suggestions on the first version of this paper.

\bibliographystyle{amsplain}
\bibliography{paperBLP}  

\end{document}